\UseRawInputEncoding
\documentclass[preprint,dvips,numbers,sort&compress]{elsarticle}
\usepackage{hyperref}
\usepackage{amsmath}
\usepackage{amssymb}
\usepackage{amsfonts}
\usepackage{amsthm}
\usepackage{arydshln}
\usepackage[toc, page]{appendix}
\usepackage{vector}
\usepackage{graphicx}
\usepackage{subfigure}
\usepackage{color}
\usepackage{paralist}
\usepackage{multirow}
\usepackage{booktabs}
\usepackage{setspace}
\graphicspath{{fingon_eps/}}

\numberwithin{equation}{section}
\def \d {\partial}
\newcommand{\defeq}{\ensuremath{\buildrel {\text{def}}\over{=}}} 
\begin{document}
%
\title{ 
Simulation of Wave in
Hypo-Elastic-Plastic Solids  
Modeled by 
Eulerian Conservation Laws}
\author[LYaddress]{Lixiang Yang\corref{c1}}
\author[RLaddress]{Robert L Lowe}
\address[LYaddress]{ 
Department of Mechanical and Aerospace Engineering, The Ohio State University, Columbus, OH, 43210
}
\address[RLaddress]{Department of Mechanical and Aerospace Engineering, University of Dayton, Dayton, OH, 45469}
\cortext[c1]{Corresponding author\\E-mail address: yang.1130@buckeyemail.osu.edu}

\begin{abstract} \noindent
This paper reports a theoretical and numerical framework to model nonlinear
waves in elastic-plastic solids.  Formulated in the Eulerian frame,  the
governing equations employed include the continuity equation, the momentum
equation, and an elastic-plastic constitutive relation. The complete governing equations are a set of
first-order, fully coupled  partial differential equations with source terms.
The primary unknowns are velocities and deviatoric stresses. By casting the governing equations into a
vector-matrix form,  we derive the eigenvalues of the Jacobian matrix to show the wave speeds.  The eigenvalues are also used to
calculate the Courant number for numerical stability.  The model equations are
solved using the Space-Time Conservation Element and Solution Element (CESE) method.  The
approach is validated by comparing our numerical results to an analytical solution for the special case of longitudinal wave motion.  

\end{abstract}
\maketitle
\section{Introduction}
\label{s:intro} \noindent
General theories of elastic-plastic waves has been treated extensively
and reviewed from various points of view 
in the literature, e.g.,
Nowacki \citep{nowacki_stress_1978},  
Clifton \citep{clifton_stress_1985}, 
Cristescu \citep{cristescu_dynamic_2007},
Craggs \citep{craggs_propagation_1957}
Herrmann \citep{oden_computational_1975}.
Numerical simulations of wave motion in elastic-plastic media
were reported by  
Buchar et al.
\citep{buchar_elastic-plastic_1983} with aid of finite element analysis,
%
Trangenstein and Collella \citep{trangenstein_higherorder_1991} studied finite deformation in elastic-plastic solids by  using hyper-elastic constitutive and kinematic evolution equation. By using higher-order Godunov scheme, they mainly studied wave propagation in hyper-elastic solids for one dimensional simulation.
Miller and Collella \citep{miller_high-order_2001} extended this work in hyper-elasticity and visco-plasticity to multiple dimension. 
Hill et al.\citep{hill_eulerian_2010} used a hybrid of the weighted essentially non-oscillatory schemes combined with explicit centered difference to solve the equations of motion expressed in an Eulerian formulation. This formulation allows for a wide range of constitutive relations.  
Giese \citep{Giese_high-resolution_2004} studied elastic-plastic wave in three space dimensions. Since the governing equations are composed of two part, he solved the flux equation by using method of transport and integrated the stress-strain relationship in time with a high order ODE solver.
In order to understand the formation of the plastic zone at the crack tip, Lin and Ballmann \citep{lin_numerical_1994}  used a characteristic-based difference method to simulate elastic-plastic wave propagation in a two-dimensional anisotropic plane strain domain. Their studies focus on small plastic deformation problems.
Tran and Udaykumar \citep{tran_particle-level_2004} developed an Eulerian, sharp interface, Cartesian grid method to simulate impact and denotation. Since energy equations are considered, the Mie-Gruneisen  equation of state is used to obtain pressure. The Essentially non-oscillatory scheme is employed to capture shocks and sharp immersed boundaries are captured by using a hybrid particle level set technique. Wang et al.\citep{wang_improved_2009} used an improved CESE method to model the impact problems of multi-material elastic plastic flows. Eulerian governing equations are adopted. Projectile nose and tail velocity are compared with experimental results. It confirms the high accuracy of CESE scheme. 
Sambasivan et al.\citep{sambasivan_simulation_2013} built a sharp interface Cartesian grid-based code to solve impact and collision problems in elasto-plastic solid medium. The ghost states are introduced to treat interface of material-to-material, material-to-void and material-to-rigid surface. Johnson-Cook material model is used in the computation.
By using a 10th order compact finite difference scheme for spatial discretization and a 4th order Runge-Kutta time marching method, Ghaisas \citep{ghaisas_unified_2018} developed a high-order Eulerian method to simulate elastic-plastic deformation as well as fluid flow. In their governing equations, derivative of the  inverse deformation gradient tensor is written to become a hyperbolic differential equation with source term. 
A multi-medium Riemann solver was proposed by Li et al.\citep{li_robust_2019} to study impact dynamics between solid and fluid. Hydro-elastoplastic constitutive material model and Mie-Gruneisen equation of state are used to close the governing equations.
Recently, Cheng et al.\citep{cheng_second-order_2020} developed a cell-centered Lagrangian scheme to model elastic-plastic flow in two-dimensional medium. Detail is given on how to construct Riemann solver for contact and elasto-plastic wave.
Since wave propagation in solid medium can be mathematically described as a set of coupled first order hyperbolic partial differential equations, Bonet et al.\citep{bonet_first_2021} casted linear momentum, the deformation gradient tensor and its co-factors into a first order partial differential equation system in Lagrangian frame. Detail mathematical properties such as hyperbolity, stability, and convergence are studied. Boscheri et al. \citep{boscheri_structure-preserving_2021} wrote deformation gradient and thermal impulse density into hyperbolic form and set the total energy to be a general summation of fluid and solid.     
Atta et al.\citep{atta_shifted_2021} suggested shifted Chebyshev polynomials of the fifth-kind as basis functions to get their approximate solutions. 
Classical large deformation elastic-plastic problems were also attempted by using virtual element method \citep{cihan_3d_2021}.
Large deformation elastic-plastic dynamics were also solved by diffuse-interface methods which were originally used to solve multiphase fluid flows\citep{jain_assessment_2023}. In the governing equations, two separate transport equations for elastic and plastic deformation tensors are combined with other conservation laws. 
Liu et al.\citep{liu_exact_2021} modeled one-dimensional multi-material elastic-plastic flow using a Riemann solver.
By using finite volume method with stationary grids in Eulerian frame, conservation laws and modified hypoelastic Wilkins model were solved to find the cause of adiabatic shear bands formation \citep{muratov_finite_2021}.
To understand seismic wave propagation, Sripanich et al.\citep{sripanich_stress-dependent_2021} studied third-order elasticity. They stated that a consistent description of conservation laws and right choice of constitutive relationship as well as elastic moduli need be used in practical scenarios.
Xiong et al.\citep{xiong_energy_2021} investigated elastic-plastic energy storage and dissipation of crystals under different strain rate impact by molecular dynamics simulations.
Nonlinear elastic-plastic wave propagation in polymers due to the action of intense energy flows were studied by Boykov et al.\citep{boykov_coupled_2022}. Finite difference method is used to solve Lagrangian description of conservation laws and hypo-elastic-plastic material relationship which is taken from the old literature.
Recently, numerical modeling of wave propagation in one dimensional elastic-plastic medium was investigated by using cell-centered Lagrangian scheme \citep{chen_cell-centered_2022}. Four types of elastic-plastic problems such as impact, tensile, piston-like and Wilkins' problems were studied.
By writing conservation of linear momentum and three geometric conservation laws (the deformation gradient, its cofactor and its determinant) into first order hyperbolic form, de Campos et al.\citep{de_campos_new_2022} used a new Updated Reference Lagrangian Smooth Particle Hydrodynamics algorithm to analyze three-dimensional large deformation elasto-plasticity problems.
In order to investigate shear band propagation, Eremin et al.\citep{eremin_microstructure-based_2022} used micro-structure based finite-difference method to model plastic flow in low carbon steel. Heuze and Stainier\citep{heuze_variational_2022} created a variational framework which will let hyperbolic conservation laws work with thermo-hyperelastic-viscoplastic constitutive equations \citep{yang_theoretical_2020}.
Aided with a high-order and pressure oscillating-free scheme, Eulerian description of conservation laws as well as a hyper-elastic framework were also used to understanding multi-material elastic-plastic flow \citep{li_diffuse-interface_2022, yang_note_2018}. Similar elastic-plastic wave propagation problems were also investigated with slightly different form of conservation laws, transport equation, or constitutive models in Eulerian frame \citep{li_HLLC-type_2022, li_complete_2022, surov_towards_2022}.
Based on multi-material diffuse interface method, Wallis et al.\citep{wallis_unified_2022} studied elasto-plastic-rigid body interaction using Eulerian conservation laws for solid-fluid interaction which were derived by Barton \citep{barton_interface-capturing_2019}.
In order to simulate large deformation and penetration problems of elasto-plastic solids, Yeom \citep{yeom_numerical_2022} numerically solved Eulerian multi-material and multi-phase flow conservation laws using a high resolution computational fluid dynamics technique on Cartesian grids.
Without considering large plastic deformation and discontinuity between multiple materials or interface between solid and fluid, the elasto-plastic constitutive model can be combined with finite element analysis to study all kinds of solid mechanics problems \citep{wu_simulation_2022, tippner_elasto-plastic_2022}.

In this paper, we will employ the Conservation Element and Solution Element (CESE) method \citep{chang_method_1995, yu_numerical_2010}, an explicit space-time finite-volume scheme, to solve a system of nonlinear elasto-dynamic model equations. As a special finite-volume method, CESE method has been formerly used to solve dynamics and combustion problems, including detonations, cavitations, flows with complex shock structures \citep{wang_direct_2005, zhang_solving_2006, jiang_spacetime_2020}. Recently, Yang \citep{yang_numerical_2021} used updated CESE scheme to study hypervelocity asteroid impacts based on an elasto-plastic flow model. During the past thirty years,many different finite-volume methods were continually applied to study solid mechanics problems \citep{cardiff_thirty_2021, cameron_discontinuities_2022}. So did the CESE method. it has been employed to solve many dynamical and vibration problems in solid structures, for instance, \citep{yu_first-order_2010, yang_numerical_2010, yang_modeling_2011, yang_velocity-stress_2011, chen_simulations_2011, lowe_modal_2014, lowe_eulerian_2016, yang_viscoelasticity_2013}. 
Since the application we are interested in is ultrasonic welding process, thermo-effect is ignored which was demonstrated by experiments. 
In this paper, we will report a novel theoretical and numerical approach to model
elastic-plastic wave motion in solids without energy equation included.  We extend the isothermal model for
modeling stress wave propagation in elastic-plastic media by using a suitable
constitutive equation to model the dynamical response of elastic-plastic
materials.

The rest of the present paper is organized as follows.  Section \ref{s:consti}
illustrates the basic formulation of constitutive relation for plasticity of
solids.  Section \ref{s:modele} summarizes the model equations, including the
continuity, momentum, and constitutive relations.  The model equations are cast
into a vector-matrix form.  The Jacobian matrix of the  one-dimensional
equations are analyzed to show the eigenvalues, which represent the wave
speeds. Linearizion of nonlinear wave equation will be given in section \ref{s:Linear}. Numerical method will be talked in section \ref{s:RRM} and section \ref{s:cese}. In section \ref{s:NR}, a numerical example will be validated by experimental results. Nonlinearity and unloading profiles of the elastic-plastic wave are further studied by our numerical method.  We then offer the limitations and concluding remarks, followed by a list of cited
references.   

\section{Hypo-Elastic-Plastic Solids} 
\label{s:consti} 
\noindent
A constitutive equation for hypo-elastic-plastic media is developed in this
section.  We first develop the constitutive equation based on
the infinitesimal theory, then generalize this model to accommodate finite
deformations.  The medium of interest is assumed
isotropic, homogeneous, non-porous, and metallic. We also adopt the customary assumptions of incompressibility of the
plastic strain, yield insensitivity to the hydrostatic part of the
stress, and, for the sake of simplicity, strain-rate-independent response.  

Based on experimental observation, infinitesimal plasticity admits the additive decomposition
\begin{equation*} 
d \varepsilon_{ij} \, = \: d \varepsilon^e_{ij} \, + \, d \varepsilon^p_{ij},  
\label{eq:plas:add_dec} 
\end{equation*} 
where $d \varepsilon^e_{ij}$ and $d \varepsilon^p_{ij}$ are the infinitesimal elastic and plastic strain increments, respectively. The elastic strain increment can be obtained from isotropic linear elasticity: 
\begin{equation*} 
d \varepsilon^e_{ij} \: = \: \frac{1 + \nu}{E} \: d \sigma_{ij} \, - \, \frac{\nu}{E} \: d \sigma_{kk} \:\! \delta_{ij},  
\label{eq:plas:ela_str_inc} 
\end{equation*}
where $d \sigma_{ij}$ is the stress increment, $\delta_{ij}$ is the Kronecker delta, $E$ is Young's modulus, and $\nu$ is Poisson's ratio.

To determine the plastic strain increment $d \varepsilon^p_{ij}$, we consider the equation of the yield surface of a material that undergoes strain-rate-independent isotropic strain-hardening: 
\begin{equation} 
F(S_{ij}, \bar{\varepsilon}^{\, p}) = 0,
\label{eq:plas:yiesur} 
\end{equation} 
where $F(S_{ij}, \bar{\varepsilon}^{\, p})$ is a scalar-valued yield function whose form is made explicit by the yield criterion, e.g., von Mises or Tresca. $S_{ij}$ is the deviatoric part of stress, i.e.,
\begin{equation*}
S_{ij} \, = \, \sigma_{ij} \, - \, \frac{1}{3} \, \sigma_{kk} \;\! \delta_{ij}.
\end{equation*}
That the yield function $F$ depends on only the deviatoric part of the stress reflects the customary assumption of yield insensitivity to hydrostatic pressure. Note that the yield surface is the union of all points in deviatoric stress space that satisfy Eq.~(\ref{eq:plas:yiesur}). The effective plastic strain $\bar{\varepsilon}^{\, p}$, defined here as
\begin{equation*} 
\bar{\varepsilon}^{\, p} = \int d \bar{\varepsilon}^{\, p},  \qquad \enspace
d \bar{\varepsilon}^{\, p} \, \defeq \, \sqrt{\frac{2}{3} \: d \varepsilon^p_{ij} \, d \varepsilon^p_{ij}} \;\! , 
\end{equation*} 
quantifies the plastic strain accumulated during the deformation history. It represents the scalar hardening parameter in (\ref{eq:plas:yiesur}).

The associated flow rule
\begin{equation} 
d \varepsilon^p_{ij} \, = \: d \lambda \: \frac{\d F}{\d S_{ij}}
\label{eq:plas:flow_rule}
\end{equation} 
implies normality of the plastic strain increment to the yield surface defined in deviatoric stress space. In Eq.~(\ref{eq:plas:flow_rule}), $d \lambda$ is a scalar function that, loosely speaking, represents the magnitude of the plastic strain increment. The loading criteria are
\begin{align} 
&F < 0  &  &\hspace{-.5in}  &  &\hspace{0.in}\hbox{elastic deformation}  \nonumber \\ 
&F = \, 0,  &  &\hspace{0.0in}\frac{\d F}{\d S_{ij}} \: d S_{ij} \, > \, 0  &  &\hspace{0.0in}\hbox{plastic loading}  \nonumber \\ 
&F = \, 0,  &  &\hspace{0.0in}\frac{\d F}{\d S_{ij}} \: d S_{ij} \, = \, 0  &  &\hspace{0.in}\hbox{neutral loading}  \nonumber \\ 
&F = \, 0,  &  &\hspace{0.0in}\frac{\d F}{\d S_{ij}} \: d S_{ij} \, < \, 0  &  &\hspace{0.in}\hbox{elastic unloading}  \nonumber 
\end{align} 
Plastic strain only accrues during plastic loading; otherwise, the plastic strain increment $d \varepsilon^p_{ij}$ vanishes. Hence, our adoption of the plastic loading criteria is tacit throughout the remainder of this section.

As strain-rate-insensitive materials harden during plastic deformation, points on the original yield surface remain on all subsequent yield surfaces.  This observation, together with a first-order Taylor series expansion of $F(S_{ij}, \bar{\varepsilon}^{\, p}),$ imply the consistency condition
\begin{equation} 
dF \, = \: \frac{\d F}{\d S_{ij}} \, d S_{ij} \, + \, \frac{\d F}{\d \bar{\varepsilon}^{\, p}} \, d \bar{\varepsilon}^{\, p} \, = \, 0.  
\label{eq:plas:consist}
\end{equation} 
Use of Eq.~(\ref{eq:plas:flow_rule}) in (\ref{eq:plas:consist}) leads to
\begin{equation}
d \lambda = - \frac{\frac{\d F}{\d S_{ij}} \, d S_{ij}}
{\frac{\d F}{\d \bar{\varepsilon}^{\, p}} \left( \frac{2}{3} \frac{\d F}{\d S_{ij}} \frac{\d F}{\d S_{ij}} \right)^{\frac{1}{2}}} \;\! .
\label{eq:plas:prop_fcn} 
\end{equation} 
We employ the von Mises yield criterion, i.e., $J_2$-flow theory, which makes the yield function $F(S_{ij},\bar{\varepsilon}^{\, p})$ in Eq. (\ref{eq:plas:yiesur}) explicit:
\begin{equation} 
F(S_{ij},\bar{\varepsilon}^{\, p}) \, = \, \frac{1}{2} \;\! S_{ij} S_{ij} \, - \, \frac{1}{3} \:\! \big[ \sigma^y (\bar{\varepsilon}^{\, p}) \big]^2 \, = \, 0,
\label{eq:plas:J2ys} 
\end{equation} 
where $J_2 = S_{ij} S_{ij}/2$ is the second invariant of the deviatoric stress (related to the energy of distortion), and $\sigma^y(\bar{\varepsilon}^{\, p})$ is the yield stress in uniaxial tension, which evolves with effective plastic strain as the material hardens during plastic deformation. For choice (\ref{eq:plas:J2ys}), it follows that
\begin{equation*} 
\frac{\d F}{\d S_{ij}} = S_{ij}, \qquad \enspace
S_{ij} \;\! d S_{ij} = \frac{2}{3} \, \sigma^y d \sigma^y, \qquad \enspace
\frac{\d F}{\d \bar{\varepsilon}^{\, p}} = - \frac{2}{3} \, \sigma^y \:\! \frac{d \sigma^y}{d \bar{\varepsilon}^{\, p}},
\end{equation*} 
and use of these results in Eq. (\ref{eq:plas:prop_fcn}) leads to
\begin{equation*} 
d \varepsilon^p_{ij} \, = \, \frac{3}{2} \frac{d \sigma^y}{\sigma^y \, \frac{d \sigma^y}{d \bar{\varepsilon}^{\, p}}} \, S_{ij}.
\end{equation*} 
For a linear strain-hardening material, the tensile yield stress increases linearly with the effective plastic strain, i.e., 
\begin{equation*} 
\sigma^y(\bar{\varepsilon}^{\, p}) \, = \, \sigma^y_o \, + \, B_{\scriptscriptstyle SH} \, \bar{\varepsilon}^{\, p},  
\end{equation*} 
where the initial tensile yield stress $\sigma^y_o$ and the strength coefficient $B_{\scriptscriptstyle SH}$ are material-dependent constants. It follows that the plastic strain increment is
\begin{equation*} 
d \varepsilon^p_{ij} \, = \, \frac{3}{2} \, \frac{d \bar{\sigma}}{B_{\scriptscriptstyle SH} \, \bar{\sigma}} \, S_{ij},  
\end{equation*}
where we have introduced the effective stress
\begin{equation*}
\bar{\sigma} \defeq \sqrt{\frac{3}{2} \, S_{ij} \;\! S_{ij}} \;\! .
\end{equation*}
Thus, the total strain increment (elastic + plastic) is
\begin{equation} 
d \varepsilon_{ij} \, = \, \frac{1 + \nu}{E} \, d \sigma_{ij} \, - \, \frac{\nu}{E} \, d \sigma_{kk} \:\! \delta_{ij} \, + \, \frac{3}{2} \, \frac{d \bar{\sigma}}{B_{\scriptscriptstyle SH} \;\! \bar{\sigma}} \, S_{ij}. 
\label{eq:plas:tot_str_inc} 
\end{equation} 
Equation (\ref{eq:plas:tot_str_inc}) can be inverted to give the deviatoric stress increment, or, equivalently, its rate 
\begin{equation} 
\dot{S}_{ij} \: = \:  2 \:\! \mu \;\! \dot{\varepsilon}_{ij} \: - \: \frac{2}{3} \, \mu \;\! \dot{\varepsilon}_{kk} \:\! \delta_{ij} \: - \: 3 \mu \, \frac{S_{kl} \;\! \dot{\varepsilon}_{kl}}{\left( \frac{B_{SH}}{2 \mu} + \frac{3}{2} \right) S_{mn} \:\! S_{mn}} \, S_{ij},  
\label{eq:plas:dev_str_inc}
\end{equation} 
where $\mu$ is the shear modulus. The infinitesimal elastic-plastic constitutive equation (\ref{eq:plas:dev_str_inc}) is generalized to the finite-deformation regime by 
\begin{inparaenum}[(i)] 
\item replacing the stress tensor $\sigma_{ij}$ with its finite Eulerian analog $T_{ij}$, the Cauchy stress,
\item replacing the infinitesimal strain increment $d \varepsilon_{ij}$ with the rate of deformation $D_{ij}$, which is the work conjugate of the Cauchy stress, and 
\item employing an objective rate $D/Dt$ of the stress.
\end{inparaenum} 
The objective rate ensures that the constitutive equation is invariant under an arbitrary superposed rigid body motion. The resulting constitutive equation is 
\begin{equation*} 
\frac{D}{Dt} \, {S}_{ij} \: = \: 2 \:\! \mu \:\! D_{ij} \, - \, \frac{2}{3} \;\! \mu \;\! D_{kk} \:\! \delta_{ij} \, - \, \beta(s) \:\! S_{kl} \:\! D_{kl} \:\! S_{ij},
\end{equation*} 
where  $s = S_{mn}S_{mn}$ and 
\begin{equation} 
\beta(s) \, = \, \left\{ 
\begin{array}{ccccc} 
0 & \hbox{if} & F < 0  \\[6pt] 
0 & \hbox{if} & F = 0 & \hbox{and} & \frac{\d F}{\d S_{ij}} \, d S_{ij} < 0  \\[6pt]
0 & \hbox{if} & F = 0 & \hbox{and} & \frac{\d F}{\d \sigma_{ij}} d \sigma_{ij} = 0  \\[6pt]
\frac{6\mu^2}{\left(3\mu + B_{SH} \right) s} & \hbox{if} & F = 0 & \hbox{and} & \frac{\d F}{\d \sigma_{ij}} d \sigma_{ij} > 0 
\end{array} 
\right.
\label{eq:plas:pl_coeff} 
\end{equation} 
Note that the four rows of Eq. (\ref{eq:plas:pl_coeff}) correspond to elastic deformation, elastic unloading, neutral loading, and plastic loading, respectively. 
 One common choice is the
Jaumann rate:
\begin{equation} 
\frac{D S_{ij}}{Dt} \; = \; \frac{\d S_{ij}}{\d t} \, + \, v_k
\frac{\d S_{ij}}{\d x_k} \, - \, W_{ik} S_{kj} \, + \, S_{ik} W_{kj},
\label{Jaumann} 
\end{equation} 
where 
$W_{ij} = (\d v_i/\d x_j - \d v_j/\d x_i)$ 
is the skew part of the velocity gradient.
The resulting constitutive equation is 
\begin{equation} 
\frac{D}{Dt}{S}_{ij} = 2 \mu D_{ij} - 
\frac{2}{3} \mu D_{kk} \delta_{ij} \; - \; 3 \mu \, \frac{S_{kl} D_{kl}}{\left(
\frac{B_{SH}}{2 \mu} + \frac{3}{2} \right) S_{mn} S_{mn}} \, S_{ij},
\label{ep_finite_1.1} 
\end{equation} 
where 
$D_{ij} = 1/2 \left( L_{ij} + L_{ji} \right)$ 
is the symmetric part of the Eulerian velocity gradient
$L_{ij} = {\d v_i}/{\d x_j}$, 
and $v_i$ is the velocity.
The first two terms in Eq.~(\ref{ep_finite_1.1}) reflect elastic contributions to the deviatoric stress, while the final term represents
the plastic contribution.
%

\section{Governing Equations}
\label{s:modele}
Based on the elastic-plastic constitutive relation, Eq. (\ref{ep_finite_1.1}), the three-dimensional governing
equations for elastic-plastic wave motion formulated in the Eulerian frame are

\vspace{0.1in} \noindent {\em Conservation of mass:}
\begin{equation} 
\frac{\d \rho}{\d t} \, + \, \frac{\d}{\d x_i} (\rho v_i) \; =
\; 0,  \label{mass} 
\end{equation} 
\noindent {\em Conservation of linear momentum:}
\begin{equation} 
\frac{\d}{\d t} \left( \rho v_i \right) \, + \, \frac{\d}{\d
x_j} \left( \rho v_i v_j + p - S_{ij} \right) \; = \; 0  \label{lin_mom}
\end{equation} 
and Cauchy stress components and pressure are $T_{ij} = - p + S_{ij}$, $p = -\sum_{i=1}^3 T_{ii}/3$.
{\em Elastic-plastic constitutive relation:} 
\begin{equation} 
\frac{D}{Dt}{S}_{ij} =  2 \mu D_{ij} - 
\frac{2}{3} \mu D_{kk} \delta_{ij} - 
\beta(s)  S_{kl} D_{kl} S_{ij},
\label{ep_finite_1} 
\end{equation} 
%
%
where  $s = S_{mn}S_{mn}$ and 
\begin{equation} 
\beta(s) \, =
\, \left\{ \begin{array}{ccccc} 
0 & \hbox{if} & F < 0 \\ 
0 & \hbox{if} & F = 0
& 
\hbox{and} & \frac{\d F}{\d \sigma_{ij}} d \sigma_{ij} < 0 
\\
0 & \hbox{if} &
F = 0 
& 
\hbox{and} & \frac{\d F}{\d \sigma_{ij}} d \sigma_{ij} = 0 
\\ 
\frac{6\mu^2}{\left(3\mu + B_{SH} \right) s} 
& 
\hbox{if} & F = 0 & \hbox{and} &
\frac{\d F}{\d \sigma_{ij}} d \sigma_{ij} > 0 \end{array} \right.
\label{pl_coeff2} 
\end{equation} 
The four rows of Eq. (\ref{pl_coeff2})
represent elastic deformation, elastic unloading, neutral loading, and
plastic loading, respectively. 
In Eqs. (\ref{mass})-(\ref{pl_coeff2}), 
$\rho$ is the density, 
$v_i$ is the velocity, 
 and 
$D_{ij} = \left( \d v_i/\d x_j + \d v_j/\d x_i \right)/2$ 
is the symmetric part of the velocity gradient,
all functions of time $t$ and the position ${\bf x} = \left\{ x_1,x_2,x_3
\right\}$ of a typical material particle in the current configuration referred
to a fixed Cartesian basis. The shear modulus $\mu$ is a prescribed material
constant. If $J_2$ flow plastic theory is adopted,  we have
\begin{align}
\frac{\d F}{\d \sigma_{ij}} d\sigma_{ij} = d (\frac{1}{2}S_{ij} S_{ij} - k^2), \label{pl_numc}
\end{align}
where $k$ is a constant for perfect plastic materials and in general it is
\begin{align*}
k^2 = \frac{1}{3} (\bar{\sigma} (\varepsilon^p_{ij}))^2. 
\end{align*}
Then Eq.(\ref{ep_finite_1}) and Eq.(\ref{pl_coeff2}) can be written as 

\begin{equation} 
\frac{D}{Dt}{S}_{ij} =  2 \mu D_{ij} - 
\frac{2}{3} \mu D_{kk} \delta_{ij} - 
\beta(s)  S_{kl} D_{kl} S_{ij},
\label{ep_finite_3} 
\end{equation} 
where  $s = S_{mn}S_{mn}$ and 
\begin{equation} 
\beta(s) \, =
\, \left\{ \begin{array}{ccccc} 
0 & \hbox{if} & F < 0 \\ 
0 & \hbox{if} &
F = 0 
& \hbox{and} &
d (\frac{1}{2}S_{ij} S_{ij} - k^2) \le 0
\\ 
\frac{6\mu^2}{\left(3\mu + B_{SH} \right) s} 
& 
\hbox{if} & F = 0 & \hbox{and} &
d (\frac{1}{2}S_{ij} S_{ij} - k^2)  > 0 \end{array} \right.
\label{pl_coeff3} 
\end{equation}

\vspace{0.225in}
\noindent {\bf {1-D Formulation}}
\vspace{0.075in}

\noindent 
For one-dimensional cases (ignore $\partial/\partial x_2, \partial/\partial x_3$),
the mass and momentum equations are
\begin{align}
  &\frac{\partial\rho}{\partial t}
    + \frac{\partial\rho v_1}{\partial x_1} = 0, \label{e:cons:mass} \\
  &\frac{\partial\rho v_1}{\partial t}
    + \frac{\partial}{\partial x_1}\left(
        \rho v_1v_1 + p - S_{11}
      \right) = 0, \label{e:cons:momentum} 
\end{align}
where $p$ is mean stress and $S_{11}$ is deviatoric stress.  Moreover, for 1D strain problem, the 3-d elastic-plastic constitutive equation can be simplified to be
\begin{align}
  &\frac{\partial\rho S_{11}}{\partial t}
    + \frac{\partial\rho v_1S_{11}}{\partial x_1}
    = \frac{4}{3}\mu\rho\left( 1 - \frac{\gamma}{1+\frac{B_{SH}}{3\mu}} \right)
      \frac{\partial v_1}{\partial x_1}, \label{e:cons:plas}
\end{align}
where
\begin{equation} 
\gamma \, =
\, \left\{ \begin{array}{ccccc} 
0 & \hbox{if} & F < 0 \\ 
0 & \hbox{if} & F = 0
& 
\hbox{and} & d (\frac{1}{2}S_{ij} S_{ij} - k^2) \le 0
\\ 
1
& 
\hbox{if} & F = 0 & \hbox{and} & d (\frac{1}{2}S_{ij} S_{ij} - k^2) > 0 \end{array} \right.
\label{pl_coeff1} 
\end{equation} 
Mentioned in previous section, to close the system of equations, we employ the following equation of state to relate pressure to density:

\begin{equation}
p = k \ln \frac{\rho}{\rho_o} + p_0. \label{e:cons:eqstate}
\end{equation}

 Eqs. (\ref{e:cons:mass})- (\ref{e:cons:plas}) could generate a hyperbolic system, which could be written in vector form as:

\begin{align}
 \frac{\partial \bvec{U}}{\partial t} + \frac{\partial \bvec{E}}{\partial x_1}
 = \bvec{H} \label{e:cons:vectform1}
\end{align}
where
\begin{align*}
\bvec{U} &= (\rho, \rho v_1, \rho S_{11})^T, \\
\bvec{E} &= (\rho v_1, \rho v_1 v_1 + p - S_{11}, \rho S_{11} v_1)^T, \\
\bvec{H} &= \left(0,0,\frac{4}{3}\mu \rho (1- \beta / (1+ B_{SH}/3 \mu))
\frac{\partial v_1}{\partial x_{1}}\right)^T. 
\end{align*}
By analyzing the eigen-structure of this hyperbolic system, we directly
calculate the speed of sound in the solid with plastic deformation. Eq. (\ref{e:cons:vectform1}) can be recast as
\begin{align}
 \frac{\partial \bvec{U}}{\partial t} 
+ \bvec{A}\frac{\partial \bvec{U}}{\partial x_1}
 = \bvec{H} 
\label{e:cons:vectform}
\end{align}
where
\begin{align}
 \bvec{A} 
    = \left(\begin{array}{ccc}
      0 & 1 & 0 \\
      -v_1^2 + \frac{k}{\rho}+\frac{S_{11}}{\rho} & 2v_1 & -\frac{1}{\rho} \\
       -v_1 S_{11} & S_{11} & v_1
    \end{array}\right)   \label{e:cons:cob1}
\end{align}
and
\begin{align}
\bvec{H} = \left(0,0,\frac{4}{3}\mu \rho (1 - \beta/(1+ B_{SH}/3\mu))\frac{\d v_1}{\d x}\right)^T
\end{align}
Alternatively, we can rewrite the above hyperbolic system by using the non-conservative variables vector
\begin{align}
\widetilde{\bvec{U}} = (\rho, v_1, S_{11})^T, \label{e:cons:nonvariables}
\end{align}
which leads to the non-conservative form:
\begin{align}
 \frac{\partial \widetilde{\bvec{U}}}{\partial t} + \widetilde{\bvec{A}}
 \frac{\partial \widetilde{\bvec{U}}}{\partial x_1}
 = \widetilde{\bvec{H}} \label{e:cons:vectform2}
\end{align}
with
\begin{align}
 \widetilde{\bvec{A}}
    = \left(\begin{array}{ccc}
      v_1 & \rho & 0 \\
      \frac{k}{\rho^2} & v_1 & -\frac{1}{\rho} \\
       0 & 0 & v_1
    \end{array}\right) \notag  \label{e:jaco}
\end{align}
and 
\begin{equation}
\widetilde{H} = \left(0,0,\frac{4}{3}\mu \left(1 - \frac{\beta}{1+ B_{SH}/3\mu}\right)\frac{\d v_1}{\d x}\right)^T
\end{equation}
By moving the source term from right side to the left side, we could transform
Eq.(\ref{e:cons:vectform2}) into
\begin{align}
 \frac{\partial \widetilde{\bvec{U}}}{\partial t} + \overline{\bvec{A}}
 \frac{\partial \widetilde{\bvec{U}}}{\partial x_1}
 = 0 \label{e:cons:vectform3}
\end{align}
where
\begin{align}
 \overline{\bvec{A}}
    = \left(\begin{array}{ccc}
      v_1 & \rho & 0 \\
      \frac{k}{\rho^2} & v_1 & -\frac{1}{\rho} \\
       0 & -\frac{4}{3}\mu (1 - \frac{\beta}{1 + B_{SH}/3\mu})\frac{\d v_1}{\d x} & v_1
    \end{array}\right) \notag  \label{e:jaco3}
\end{align}
which is suitable for assessing the eigenstructure of the system of equations.
The eigenvalues of matrix $\overline{A}$ can be readily
derived and they are
\begin{equation}
\lambda_1 = v_1, \lambda_{2,3}= v_1\pm c = v_1 \pm \sqrt{\frac{k + \frac{4}{3}\mu(1 - \frac{\beta}{1+B_{SH}/3\mu})}{\rho}},
\label{e:cons:eigenvalue1}
\end{equation}
where $k$ is bulk modulus, $\mu$ is shear modulus.  It can be seen that  is
$\sqrt{(k + \frac{4}{3}\mu)/\rho}$ by letting $\beta = 0 $.  It can be seen
that this plastic wave speed is slower than elastic wave speed in the bulk
material. In particular, for the elastic-perfectly plastic materials, i.e. $\beta=1$ and
$B_{SH} = 0$, the plastic wave speed is given by
\begin{equation}
c =\sqrt{\frac{k}{\rho}}. 
\label{e:cons:eigenvalue2}
\end{equation}
In the rest of the present paper, the above formulation will be numerically
solved by the CESE method.  The computational conditions and numerical results
and remarks will be illustrated in the following sections.
\section{Linearization of nonlinear plastic wave and wave speed illustration}
\label{s:Linear}
\begin{figure}
\centering
  \includegraphics[scale = 0.3]{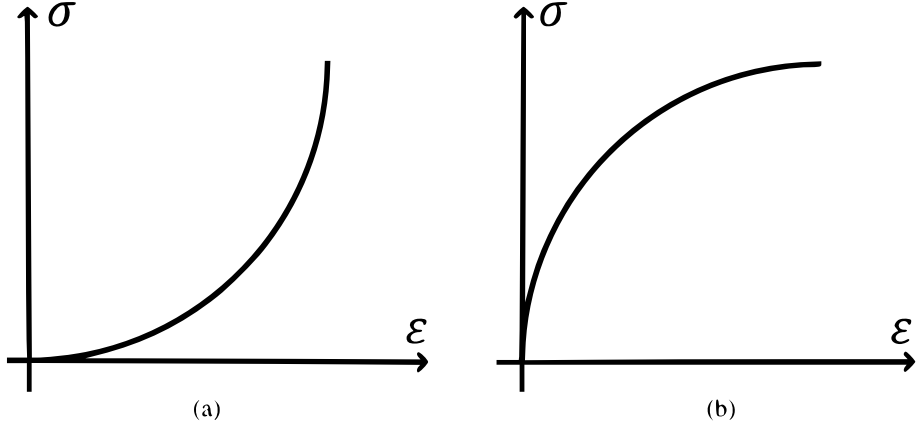}
\caption{
A schematic of one dimensional stress-strain curve. (a) A plastic stress-strain curve for a work-hardening material.
(b) A plastic stress-strain curve concave towards the stress axis.
}
\label{f:onedimplconstitutive}
\end{figure}

For the purpose of comparison to the small-amplitude numerical solution of the
nonlinear wave equations, we obtain the analytical solution of the linearized problem. We
linearize the nonlinear governing equations
(\ref{e:cons:mass}), (\ref{e:cons:momentum}) and (\ref{e:cons:plas}) as follows. We
expand $\rho(x_1,t)$, $v(x_1,t)$ and $S_{11}$in $\epsilon$ about the rest state
$\rho^0$, $v^0$ and $S_{11}^0$:
\begin{align}
\begin{aligned} 
\rho(x_1,t) &= \rho^0 + \epsilon \rho^1(x_1,t) + \epsilon^2
\rho^2(x_1,t) +\ldots,    
\\
v(x_1,t) &= v^0 + \epsilon v^1(x_1,t)
+\epsilon^2v^2(x_1,t) + \ldots, 
\\
S_{11}(x_1,t) & = S_{11}^0 + \epsilon
S_{11}^1(x_1,t) + \epsilon^2 S_{11}^2(x_1,t) + \ldots, 
\label{e:cons:linear}
\end{aligned} 
\end{align}
where $\rho^i(x_1,t)$, $v^i(x_1,t)$ and $S_{11}^i(x_1,t)$, $i = 1,2,...$ are
the $i$-th corrections of density, velocity and deviatoric stress, respectively,
and $\epsilon$ is a small, dimensionless, positive, scalar quantity. Inserting
the expansions into Eq. (\ref{e:cons:mass}), Eq.
(\ref{e:cons:momentum}) and Eq. (\ref{e:cons:plas}), selecting $\rho^0 =
\rho_0$, $v^0 = 0$ and $S_{11}^0 = 0$ as the rest state, the following linear
equations are obtained from the order-$\epsilon$ problem:
\begin{align}
\begin{aligned}
\frac{\d \rho^1}{\d t} + \rho_0 \frac{\d v^1}{\d x_1} = 0,\\
\rho_0 \frac{\d  v^1}{\d t} + \frac{K}{\rho_0} \frac{\d \rho^1}{\d x} - \frac{\d S_{11}^1}{\d x} = 0,\\
\frac{\d S_{11}^1}{\d t}  = \frac{4}{3}\mu [1- {\beta}/(1+ {B_{SH}}/{3 \mu})]\frac{\d v^1}{\d x}.
\label{e:cons:linear1}
\end{aligned}
\end{align} 
Equations (\ref{e:cons:linear1}) are combined to recover second order
linear wave equations in velocity:
\begin{align}
\frac{\d^2 v^1}{\d t^2} = c^2 \frac{\d^2 v^1}{\d x^2}.
\label{e:cons:2ndorderwave} \end{align}
where
\begin{align*} 
c = \sqrt{\frac{K + \frac{4}{3}\mu [1 - {\beta}/(1+
{B_{SH}}/{3\mu})]}{\rho_0}} 
\label{e:cons:lnwvspd} 
\end{align*}
is the plastic wave propagation speed. If we integrate
Eq. (\ref{e:cons:2ndorderwave}) with respect to time and set the arbitrary
function of $x$ which arise to be zero, we recover
\begin{align}
\frac{\d^2 u}{\d t^2} = c^2 \frac{\d^2 u}{\d x^2} 
\end{align}
with $u(x,t)$ the axial displacement component.  It can be seen in
Eq.(\ref{e:cons:lnwvspd}) that plastic wave speed depends on plastic modulus
$B_{SH}$ which is constant in metals with linear hardening plastic materials, 
a reasonable first approximation for strain hardening. So the variation of the
velocity of plastic wave propagation, as a function of strain, is governing
by the slope of the plastic stress-strain curve. For instance,
the stress-strain curve of some materials takes the form illustrated in
Fig.(\ref{f:onedimplconstitutive}a), i.e. where $\sigma d^2\sigma/d \epsilon^2 < 0$. For
such case, the plastic wave speed $c(\epsilon)$ increases when the stress increase. In
the impact case, since the stress will increase at the impact end, the waves
generated by the impact will propagate with continually increasing velocities. In this case the distance
between the wave fronts becomes short during propagation and there is a
tendency to form shock waves. Detail analysis of infinitesimal plastic wave
speed was shown by Cristescu \citep{cristescu_dynamic_2007} (see Appendix). However, for some rubbers, soils
and certain metals, the constitutive equation takes the form shown in Fig.(\ref{f:onedimplconstitutive}b), i.e. for which the slope decreases continuously
($\sigma d^2 \sigma /d \epsilon^2 > 0$) for any strain. The speed of wave
propagation will decrease when the stress increases. 
This also means that the wave fronts will become wide during their propagation.
Therefore, an expansion wave will be formed. Several numerical examples will be given in the results section. 

\section{Numerical Solution by Radial Return Mapping}
\label{s:RRM}
Based on the consistency condition, the effective stress must be constrained to
always fall either within or on the yield surface. Different from 
infinitesimal plasticity, in solving finite plasticity problem,
the typical numerical method is radial return algorithm
\citep{tran_particle-level_2004}. This algorithm includes two steps: (i) predict a trial stress 
using the elastic part of Eq. (\ref{e:cons:plas}),

\begin{equation}
\frac{\partial \rho S_{11,tr}}{\partial t} + 
\frac{\partial \rho u S_{11,tr}}{\partial x} 
= \frac{4}{3}\mu\rho\frac{\partial u}{\partial x} 
\label{e:cons:RRA1}
\end{equation}
and (ii) correct the trial stress to be true elastic-plastic stress by pulling
the trial stress back to the yield surface:
\begin{align} 
S_{11,tr} &= S_{11,tr}-\frac{S_{11,tr}}{\mid
S_{11,tr}\mid} \cdot \frac{\mid S_{11,tr} \mid - \mid S_{11,pre} \mid} {1 +
\frac{B}{3 \mu}} 
\label{e:cons:rra12} 
\end{align}
where $S_{11.tr}$ is the trial stress predicted by assuming purely elastic
deformation, $S_{11,pre}$ is the true elastic-plastic stress calculated at the
previous time step. 

\section{The CESE Method}
\label{s:cese}

Conventional finite volume methods are formulated according to a flux balance
over a fixed spatial domain. The conservation laws state that the rate of
change of the total amount of a substance contained in a fixed spatial domain,
i.e., the control volume $V$, is equal to the flux of that substance across the
boundary of $V$, denoted as $S(V)$. Consider the differential form of a
conservation law as follows: \begin{align} \frac{\partial u}{\partial t} +
\bigtriangledown \cdot \bvec{f} = 0 \label{e:cese:1} \end{align} where $u$ is
density of the conserved flow variable, $f$ is the spatial flux vector. By
applying Reynold's transport theorem to the above equation, one can obtain
the integral form as: \begin{align} \frac{\partial}{\partial t}\int_V u dV +
\oint_{S(V)}\bvec{f} \cdot d \bvec{S} = 0 \label{e:cese:2} \end{align} where
$dV$ is a spatial volume element in $V$, $d \bvec{s}=d \sigma \bvec{n}$ with
$d\sigma$ and $\bvec{n}$ being the area and the unit outward normal vector of a
surface element on $S(V)$ respectively. By integrating Eq.(\ref{e:cese:2}), one
has \begin{align} \left[\int_V u dV \right]_{t = t_f} - \left[\int_V u dV
\right]_{t = t_s} + \int^{t_f}_{t_s} \left(\oint_{S(V)}\bvec{f} \cdot d
\bvec{S}\right)dt = 0 \label{e:cese:3} \end{align} The discretization of
Eq.(\ref{e:cese:3}) is the focus of the conventional finite-volume methods. In
particular, the calculation of the flux terms in Eq.(\ref{e:cese:3}) would
introduce the upwind methods due to the nonlinearity of the convection terms in
the conservation laws. 

In the CESE method, we do not use the above formulation based on the Reynolds
transport theorem. Instead, the conservation law is formulated by treating
space and time on an equal-footing. This unified treatment of space and
time allows a consistent integration in space-time and thus
ensures local and global flux balance. This chapter briefly illustrates the
CESE method in one-spatial dimension. 

\subsection{One-Dimensional CESE Method} 
\begin{figure} \centering
\includegraphics[scale = 0.6]{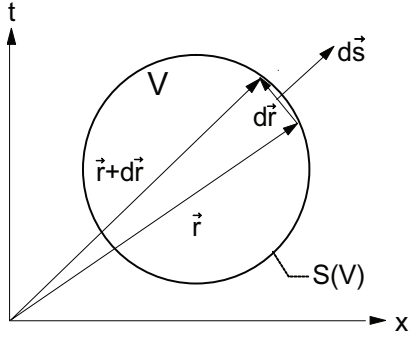} \caption{A schematic of space-time
integral of the CESE method} \label{cese-1d1} \end{figure}

To proceed, let time and space be the
two orthogonal coordinates of a space-time system, i.e.,$x_1 = x$  and $x_2 =
t$. They constitute a two-dimensional Euclidean space $E_2$. Define $h \equiv
(f,u)$, then by using the Gauss divergence theorem, Eq.(\ref{e:cese:1})
becomes \begin{align} \int_{\partial \Omega} \bvec{h} \cdot d \bvec{s} = 0
\label{e:cese:4} \end{align} Equation (\ref{e:cese:4}) states that the total
space-time flux $\bvec{h}$ leaving the space-time volume through its surface
vanishes. Refer to Figure \ref{cese-1d1} for a schematic of
Eq.(\ref{e:cese:4}). To integrate Eq. (\ref{e:cese:4}) we employ the CESE
method \citep{chang_method_1995}.   

\begin{figure}
\centering
\subfigure[]{
  \includegraphics{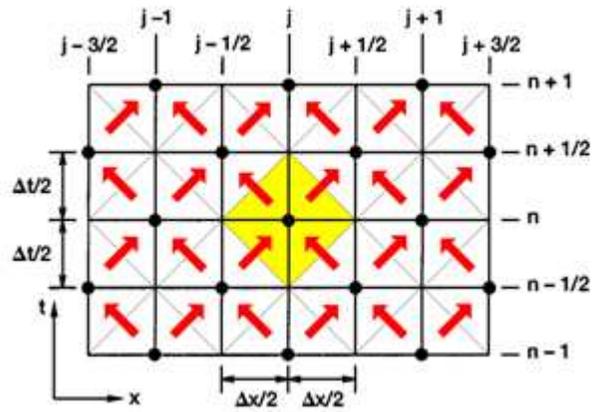}
  \label{f:cese:zigzag}
}
\subfigure[]{
  \includegraphics{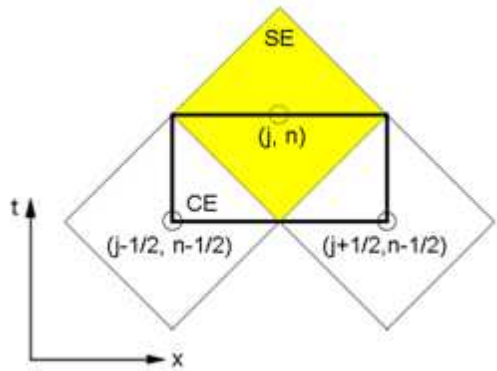}
  \label{f:cese:ceint}
}
\caption[A schematic of the CESE method in one spatial dimension.
\subref{f:cese:zigzag} Zigzagging SEs.
\subref{f:cese:ceint} Integration over a CE to solve $u_i$ and $(u_x)_i$ at the
new time level.]{
A schematic of the CESE method in one spatial dimension.
\subref{f:cese:zigzag} Zigzagging SEs.
\subref{f:cese:ceint} Integration over a CE to solve $u_i$ and $(u_x)_i$ at the
new time level.
}
\label{f:cese-1d2}
\end{figure}

In the CESE method, separated definitions of Solution Element (SE) and
Conservation Element (CE) are introduced. In each SE, solutions of unknown
variables are assumed continuous and a prescribed function is used to represent
the profile. In the present calculation, a linear distribution is used. Over
each CE, the space-time flux in the integral form, Eq.(\ref{e:cese:4}), is
imposed.  Figure \ref{f:cese-1d2} shows the space-time mesh and the
associated SEs and CEs. Solutions of variables are stored at mesh nodes which
are denoted by filled circular dots. Since a staggered mesh is used, solution
variables at neighboring SEs leapfrog each other in time-marching calculation.
The SE associate with each mesh node is a yellow rhombus. Inside the SE, the
solution variables are assumed continuous. Across the interfaces of neighboring
SEs, solution discontinuities are allowed. In this arrangement, solution
information from one SE to another propagates only in one direction, i.e.,
toward the future through the oblique interface as denoted by the red arrows.
Through this arrangement of the space-time staggered mesh, the classical Riemann
problem has been avoided. Figure \ref{f:cese:ceint} illustrates a rectangular
CE, over which the space-time flux conservation is imposed. This flux balance
provides a relation between the solutions of three mesh nodes: $(j,n)$,
$(j-1/2, n-1/2)$, and $(j+1/2, n-1/2)$. If the solutions at time step $n-1/2$
are known, the flux conservation condition would determine the solution at
$(j,n)$.      

In the present research, many differential equations have source terms. Thus,
we consider the one-dimensional equations with source terms: 
\begin{equation}
\frac{\partial u_m}{\partial t} + \frac{\partial f_m}{\partial x} = s_m,
\label{e:cesec2} \end{equation} 
where $m = 1,2,3$ and the source term $s_m$ are
functions of the unknowns $u_m$ and their spatial derivatives. For any $(x,t)
\in \mathrm{SE}(j,n)$, $u_m(x,t)$, $f_m(x,t)$ and $\mathbf{h}_m(x,t)$, are
approximated by $u^*(x,t;j,n)$, $f^*(x,t;j,n)$, and $\mathbf{h}^*(x,t;j,n)$.
By assuming linear distribution inside an SE, we have \begin{align*} &u_m^*(x,
t; j, n) = \\ &\quad (u_m)_j^n + (u_{mx})_j^n(x-x_j) + (u_{mt})_j^n(t-t^n), \\
&f_m^*(x, t; j, n) = \\ &\quad (f_m)_j^n + (f_{mx})_j^n(x-x_j) +
(f_{mt})_j^n(t-t^n), \\ &\mathbf{h}^*_m(x,t;j,n) = (f_m^*(x,t;j,n),
u_m^*(x,t;j,n)), \end{align*} where \begin{align*} (u_{mx})_j^n &=
\left(\frac{\partial u_m}{\partial x}\right)_j^n, \\ (f_{mx})_j^n &=
\left(\frac{\partial f_m}{\partial x}\right)_j^n = (f_{m,l})_j^n(u_{lx})_j^n,
\\ (u_{mt})_j^n &= \left(\frac{\partial u_m}{\partial t}\right)_j^n =
-(f_{mx})_j^n = -(f_{m,l})_j^n(u_{lx})_j^n, \\ (f_{mt})_j^n &=
\left(\frac{\partial f_m}{\partial t}\right)_j^n \\ &=
(f_{m,l})_j^n(u_{lt})_j^n = -(f_{m,l})_j^n(f_{l,p})(u_{px})_j^n, \end{align*}
and $(f_{m,l})_j^n \equiv (\partial f_m/\partial u_l)_j^n$ is the Jacobian
matrix.  Assume that, for any $(x,t) \in \mathrm{SE}(j,n)$, $u_m =
u_m^*(x,t;j,n)$ and $f_m = f_m^*(x,t;j,n)$ satisfy Eq.~(\ref{e:cesec2}), i.e.,
\begin{equation} \frac{\partial u_m^*(x,t;j,n)}{\partial t} + \frac{\partial
f_m^*(x,t;j,n)}{\partial x} = s_m^*(x,t;j,n), \label{e:cesec4} \end{equation}
where we assume that $s_m^*$ is constant within $\mathrm{SE}(j,n)$, i.e.,
$s_m^*(x,t;,j,n) = (s_m)_j^n$. Eq.~(\ref{e:cesec4}) becomes \begin{equation}
(u_{mt})_j^n = -(f_{mx})_j^n + (s_m)_j^n. \label{e:cesec5} \end{equation}
Since$(f_{mx})_j^n$ are functions of $(u_m)_j^n$ and $(u_{mx})_j^n$; and
$(s_m)_j^n$ are also functions of $(u_m)_j^n$, Eq.~(\ref{e:cesec5}) implies
that $(u_{mt})_j^n)$ are also functions of $(u_m)_j^m$ and $(u_{mx})_j^m$.
Aided by the above equations, we determine that the only unknowns are
$(u_m)_j^n$ and $(u_{mx})_j^n$ at each mesh point $(j,n)$.

Next, we impose space-time flux conservation over $\mathrm{CE}(j, n)$ to
determine the unknowns $(u_m)_j^n$.  Refer to Fig.~\ref{f:cese:ceint}.  Assume
that $u_m$ and $u_{mx}$ at mesh points $(j-1/2, n-1/2)$ and $(j+1/2, n-1/2)$
are known and their values are used to calculate $(u_m)_j^n$ and $(u_{mx})_j^n$
at the new time level $n$.  By enforcing the flux balance over $\mathrm{CE}(j,
n)$, i.e.,  
\begin{align*}
  \oint_{S(\mathrm{CE}(j,n))}\bvec{h}_m^*\cdot d\bvec{s} = 
  \int_{\mathrm{CE}(j,n)} s_m^* d\Omega,
\end{align*}
one obtains
\begin{align}
  &(u_m)_j^n - \frac{\Delta t}{4}(s_m)_j^n = \frac{1}{2}\Bigl[
      (u_m)_{j-1/2}^{n-1/2} + (u_m)_{j+1/2}^{n-1/2} \notag \\
  &\quad
    + \frac{\Delta t}{4}(s_m)_{j-1/2}^{n-1/2}
    + \frac{\Delta t}{4}(s_m)_{j+1/2}^{n-1/2} \notag \\
  &\quad
    + (p_m)_{j-1/2}^{n-1/2} - (p_m)_{j+1/2}^{n-1/2}
  \Bigl],
\label{e:cese:sol:um}
\end{align}
where
\begin{align*}
  (p_m)_j^n = \frac{\Delta x}{4}(u_{mx})_j^n
    + \frac{\Delta t}{\Delta x}(f_m)_j^n
    + \frac{\Delta t^2}{4\Delta x}(f_{mt})_j^n.
\end{align*}
Given the values of the marching variables at the mesh nodes $(j-1/2, n-1/2)$
and $(j+1/2,n-1/2)$, the right-hand side of Eq.~(\ref{e:cese:sol:um}) can be
explicitly calculated.  Since $(s_m)_j^n$ on the left hand side of
Eq.~(\ref{e:cese:sol:um}) is a function of $(u_m)_j^n$, we use Newton's method
to solve for $(u_m)_j^n$. The initial guess of the Newton iterations is
\begin{align*}
  &(\bar{u}_m)_j^n = \frac{1}{2}\Bigl[
      (u_m)_{j-1/2}^{n-1/2} + (u_m)_{j+1/2}^{n-1/2} \\
  &\quad + \frac{\Delta t}{4}(s_m)_{j-1/2}^{n-1/2}
  + \frac{\Delta t}{4}(s_m)_{j+1/2}^{n-1/2} \\
  &\quad + (p_m)_{j-1/2}^{n-1/2} - (p_m)_{j+1/2}^{n-1/2}
  \Bigl],
\end{align*}
i.e., the explicit part of the solution of $(u_m)_j^n$. 

The solution procedure for $(u_{mx})_j^n$ at node $(j,n)$ follows the standard
$a$-$\varepsilon$ scheme \citep{chang_method_1995} with $\varepsilon = 0.5$. 
To proceed, we let
\begin{align}
  (u_{mx})_j^n = \frac{(u_{mx}^+)_j^n + (u_{mx}^-)_j^n}{2},
  \label{e:cese:sol:umx}
\end{align}
where
\begin{align*}
  (u_{mx}^{\pm})_j^n &= \pm\frac{(u_m)_{j\pm1/2}^n - (u_m)_j^n}{\Delta x/2}, \\
  (u_m)_{j\pm1/2}^n &= (u_m)_{j\pm1/2}^{n-1/2}
    + \frac{\Delta t}{2}(u_{mt})_{j\pm1/2}^{n-1/2}.
\end{align*}
For solutions with discontinuities, Eq.~(\ref{e:cese:sol:umx}) is replaced by a
re-weighting procedure to add artificial damping at the jump
\begin{align*}
  (u_{mx})_j^n = W\left( (u_{mx}^-)_j^n, (u_{mx}^+)_j^n, \alpha \right),
\end{align*}
where the re-weighting function $W$ is defined as:
\begin{align*}
  W(x_-, x_+, \alpha) = \frac{|x_+|^{\alpha}x_- + |x_-|^{\alpha}x_+}
                             {|x_+|^{\alpha} + |x_-|^{\alpha}},
\end{align*}
and $\alpha$ is an adjustable constant.  The complete discussion of the
one-dimensional CESE  method can be found in \citep{chang_new_1991,
chang_method_1995}. The above method with CE and SE defined as in
Fig.~\ref{f:cese-1d2} is useful for solving the hyperbolic PDEs with non-stiff source
terms.

\section{Numerical Results}
\label{s:NR}
\subsection{Semi-infinite Domain Impact Analysis}
\begin{figure}
\centering
\includegraphics[scale = 0.3]{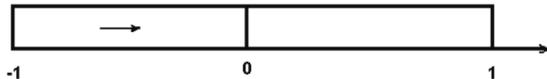}
\caption{
Initial condition of the one-dimensional impact problem
}

\label{f:cese_f1}
\end{figure}

We consider a one-dimensional copper bar with an initial speed $u= 40 m/s$
hitting a stationary copper bar.  Refer to Figure \ref{f:cese_f1}.
The initial pressures $p$ and deviatoric stress component $S_{11}$ in both
copper bars are zero. The material properties of copper are listed in Table
\ref{e:cons:table1}.

\begin{table} 
\caption{Material properties of copper} \centering
\begin{tabular}{ccccc} 
\hline
\hline 
$k(GP_a)$ &  $\rho_o (kg/m^3)$    &   $\mu
(GP_a)$  & $E (GP_a)$   & $\sigma_y (MP_a)$  
\\ 
\hline 
140   & 8930 & 45 & 122 & 90 \\ [1ex] 
\hline 
\end{tabular} 
\label{e:cons:table1} 
\end{table}
We assume the material is elastic-perfectly plastic, i.e.  the yielding stress
always equals to the initial yield stress without hardening ($B_{SH} = 0$ in Eq.
(\ref{e:cons:linear1})).  The boundary conditions at the left end of initially moving
copper bar and the right end of initially static copper bar are set as the
non-reflective boundary conditions. The focus of the present impact problem is
the interactions between the moving bar and the initially static bar. The
non-reflective boundary conditions at the two far ends allows clear observation
of wave evolution initiated  from the impact.
 
The computational domain is 2 meters, which is uniformly discretized into 400
numerical cells. The time step for the time marching calculation is 0.6 $\mu
s$. Based on the known size of spatial grid, the time increment, and the
longitudinal plane wave speed in copper bar, the CFL number in computation is
controlled to be about 0.6. The physical duration of wave propagation in
computation is 0.17 ms.

\begin{figure}
\centering
\subfigure[]{
  \includegraphics[scale = 0.4]{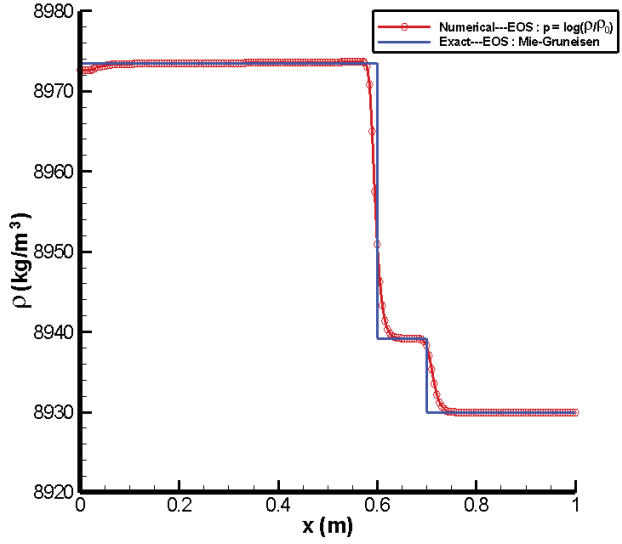}
  \label{f:cese:f2}
}
\subfigure[]{
  \includegraphics[scale = 0.4]{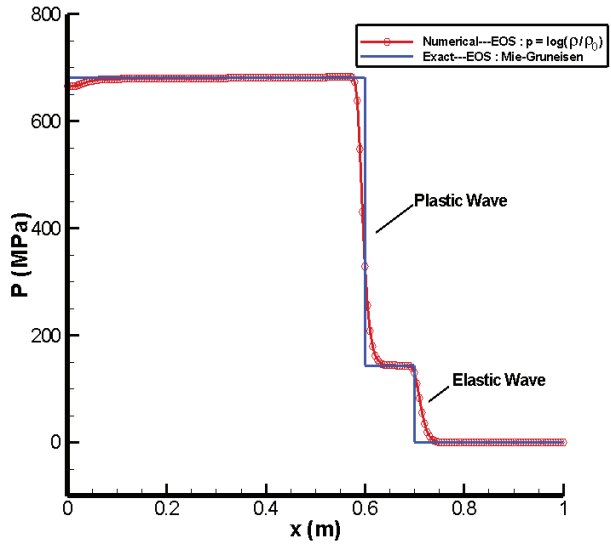}
  \label{f:cese:f3}
}

\caption{
A snapshot of density and pressure at $t = 0.17 ms$ in an initially static  copper
bar. The CESE numerical result by using the isothermal model is compared to the
exact solution by Udaykumar et al. (\citep{tran_particle-level_2004})
\subref{f:cese:f2} a snapshot of density at $t = 0.17 ms$.
\subref{f:cese:f3} a snapshot of pressure at $t = 0.17 ms$.
}

\label{f:cese2}
\end{figure}

In Fig.\ref{f:cese:f2} and Fig.\ref{f:cese:f3}, red lines with symbols represent the numerical
solutions of density and pressure by the CESE method. The blue solid lines in
these two figures represent the exact solutions by Udaykumar et
al. (\citep{tran_particle-level_2004}). They used the Mie-Gruneisen equation as the
equation of state to relate internal energy, pressure and density. 
Fig.\ref{f:cese:f2} and Fig.\ref{f:cese:f3} show that the numerical 
solutions by solving our isothermal model
equations compare well with the analytical
solution (\citep{tran_particle-level_2004}) in terms of the wave locations and
amplitude for both the plastic wave and the precursive elastic wave. This
 agreement between numerical solutions and
the analytical solution shows that in the range of low-impact velocity, the
material response simulated by the simple equation of state 
asymptotically approaches that simulated by the Mie-Gruneisen
equation.

Both the exact solution and the numerical solution show that the elastic wave
is faster than the plastic wave. This is consistent with elastic-plastic wave
speed shown in Eq.(\ref{e:cons:eigenvalue1}). In a solid with pure elastic
deformation, the wave speed is $c = \sqrt{[k + (4/3)/\mu]/\rho}$. When the
deformation involves perfect plasticity, the wave speed is $c = \sqrt{k/\rho}$,
which is lower than elastic wave speed.

To proceed, consider a one-dimensional copper bar with an initial speed $u= 30 $m/s hitting a stationary copper bar. The materials are assumed to be elastic and isotropic
linear strain hardening plastic. Three cases are considered in this example. Case 1:
the material has one effective yielding stress, e.g., 60 MPa. When strain is small, the
material is elastic, e.g., elastic Young's modulus (122 GPa), shear modulus
(45 GPa) and bulk modulus (140 GPa). After equivalent stress surpasses
effective yielding stress, the material properties will change once, 
which is decided by the effective stress-strain curve; Case 2: The
material will have two effective yielding stresses, e.g., 60 MPa and 80
MPa. If stress is larger than yielding stress, the material properties will
change twice.  Case 3: The material has 4
effective yielding stresses, e.g., 60 MPa, 80 MPa, 100 MPa, 120 MPa. These multiple
effective yielding stresses can be considered as discretized hardening stages for plastic deformation.
With the materials' properties defined, we use CESE method to calculate elastoplastic
wave propagation. In the computation, $\Delta x = 0.66 mm$ and $\Delta t = 1.3\times
10^{-7}$s. The initial pressures $p$ and deviatoric stress component $S_{11}$
in both copper bars are zero.
\begin{figure*}
\centering
\subfigure[]{
  \includegraphics[scale = 0.3]{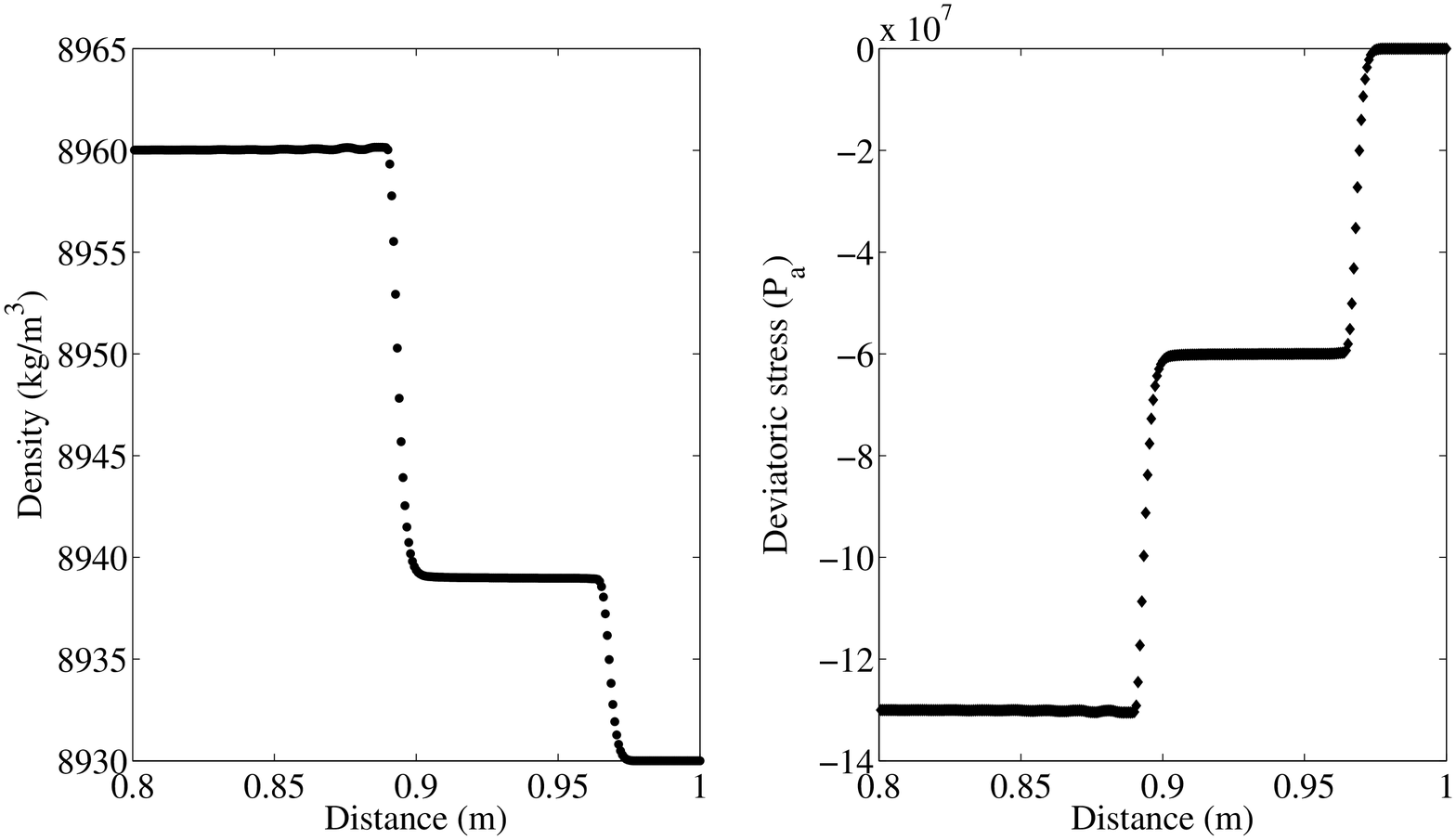}
  \label{f:plastic:1}
}
\subfigure[]{
  \includegraphics[scale = 0.3]{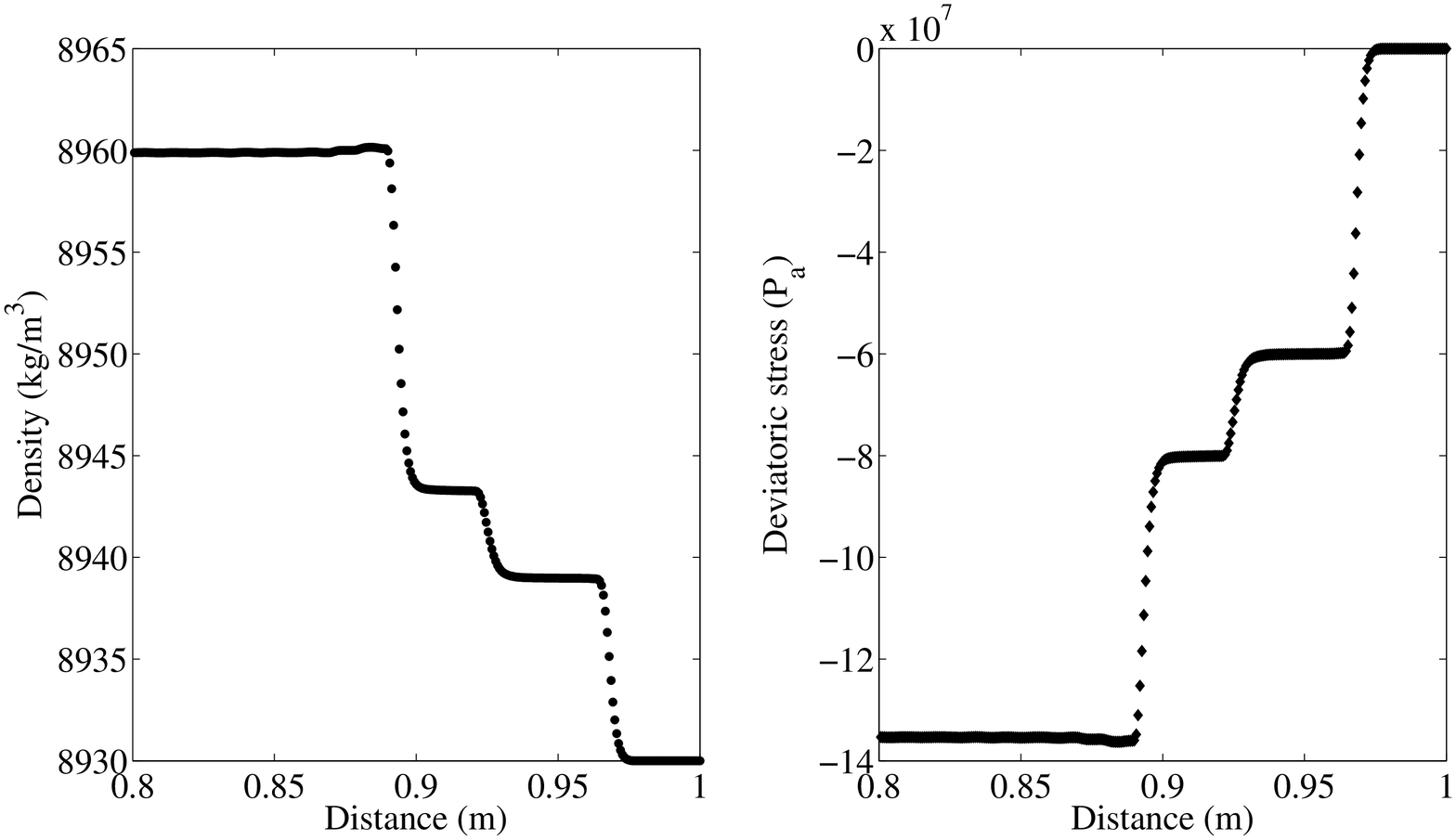}
  \label{f:plastic:2}
}
\subfigure[]{
  \includegraphics[scale = 0.3]{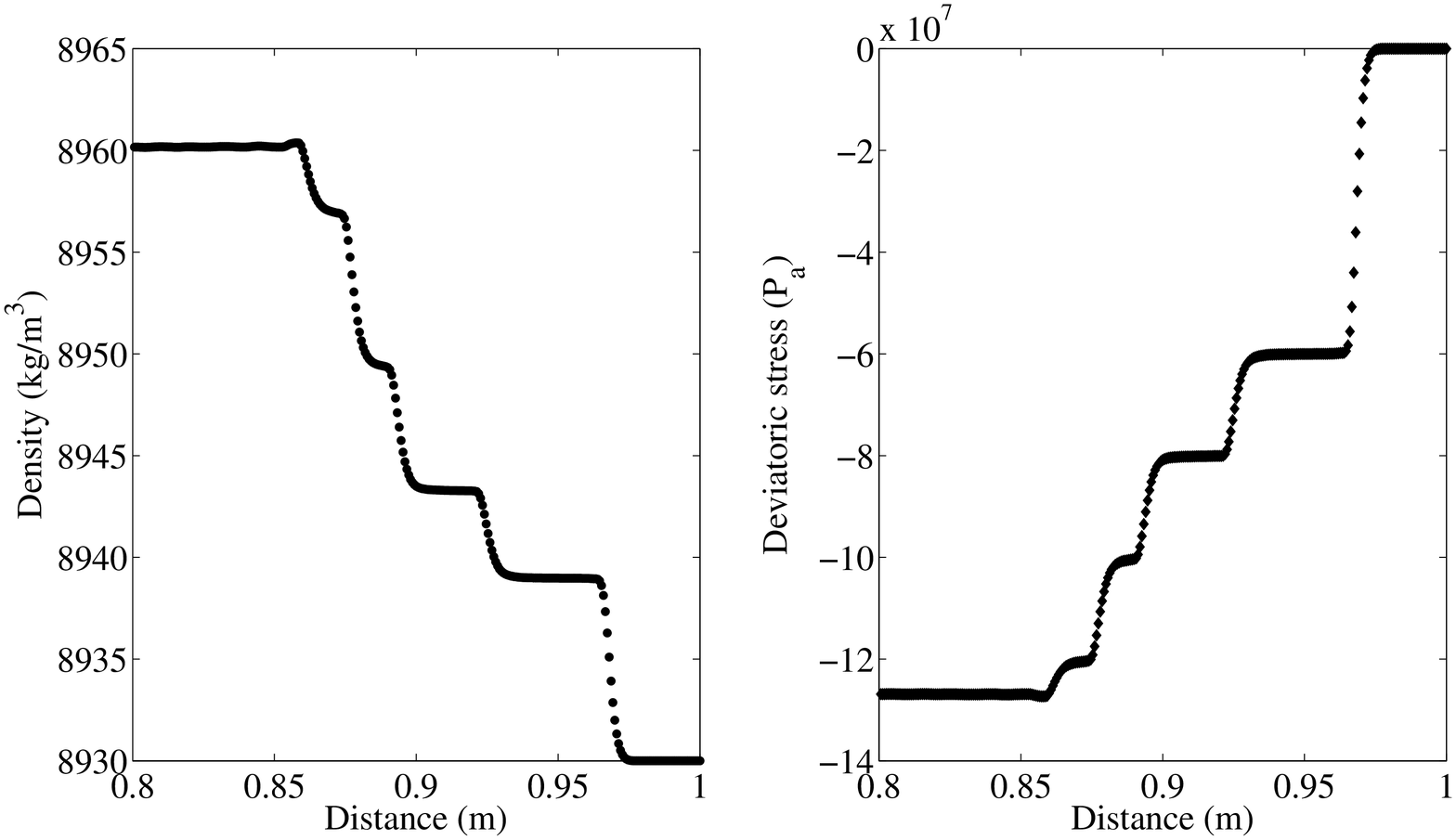}
  \label{f:plastic:3}
}
\caption{Density and stress profiles of plastic wave propagation in bulk materials for 3 cases:
(top) material has one isotropic hardening effective yielding stress; (middle)
material has two isotropic hardening effective yielding stress; (bottom)
material has four isotropic hardening effective yielding stress.}
\label{f:np:pl-5}
\end{figure*}
Fig.\ref{f:np:pl-5} shows density and deviatoric stress
profiles for three cases. Elastoplastic wave propagates to the right.  It
was demonstrated there are two wave fronts, one elastic
wave and one plastic wave, for case 1. From the discussion earlier, we know that
elastic wave front is faster than plastic wave front. In
case 2, there are 3 wave fronts, where the fastest wave is elastic wave. The other two wave fronts are plastic wave. The plastic wave speeds will depend
on the two plastic hardening slopes of effective stress and strain curve. Similarly, in case 3, there are five wave fronts, e.g., four plastic wave fronts and one elastic wave front. When the effective 
yielding stress number increases, eventually, the effective stress and strain curve 
becomes smooth. The elastoplastic wave fronts could become infinite. Therefore,
 an expansion wave will be formed. 
\begin{figure*}
\centering
  \includegraphics[scale = 0.35]{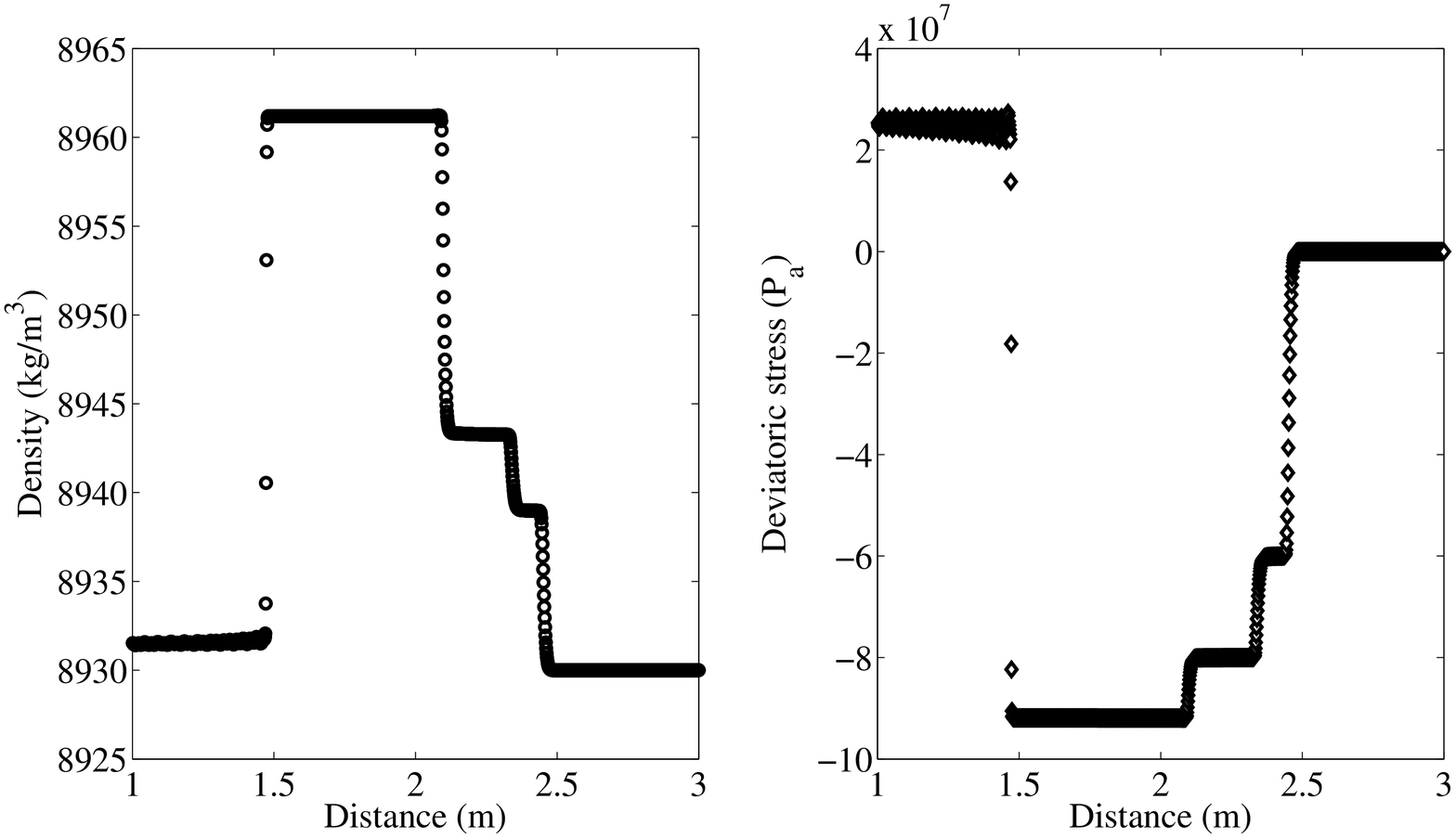}

\caption{The density profile (left) and the stress profile (right) of stress wave propagation in bulk elastic-hardening plastic materials. The unloading elastic stress wave propagates to the right with elastic wave speed. The two plastic stress wave fronts propagate to the right with two slow plastic wave speeds.}
  \label{f:plastic7}

\end{figure*}
\subsection{Unloading of Plastic Wave}
\noindent
In this numerical example, we will test the elastic and plastic wave unloading.
It is well known that materials can incur residual stress during unloading
from plastic zone in static state. In this example, the unloading of plastic 
wave propagation will be demonstrated. The numerical example is set to be the
same as the former impact analysis except that the length of the impact bar is short and finite.
The left side of the bar is stress free. Stress wave will switch from the tensive to
compressive and vice versa at the free end.

Fig.(\ref{f:plastic7}) shows a short length copper bar impacting a semi-infinite copper
bar. After contact, elastoplastic wave is generated from the contact area. 
One is left running wave, the other one is
right running wave. The left running wave will be reflected from the left end
 and join the right running wave. In Fig.(\ref{f:plastic7}), all of the elastic and plastic waves propagate rightward. The elastic wavefront is the fastest, followed by two plastic wavefronts.
 But there is only one elastic wave speed
in the unloading zone. It demonstrates that unloading path will be parallel to
the elastic curve in the constitutive relationship.  
\subsection{Ultrasonic Plastic Wave} 
\begin{figure*} 
\centering \includegraphics[scale = 0.35]{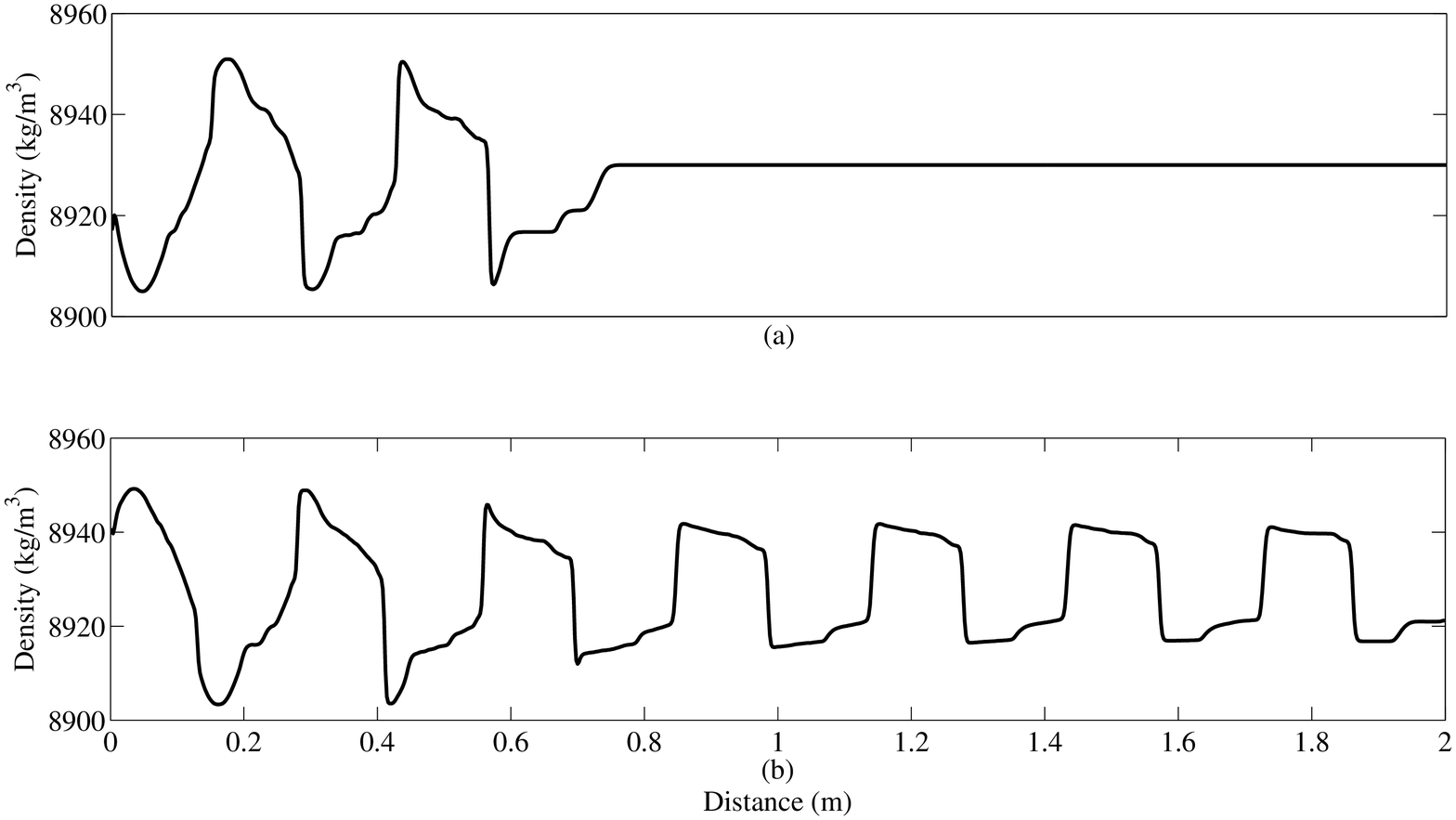}
\caption{Density profiles of stress wave propagation in bulk elastic-hardening
plastic materials. The left end is applied with a sinusoidal force in 10 kHz
frequency. Amplitude of the applied force is larger than yielding stress of
copper.} 
\label{f:plastic8} 
\end{figure*}

\begin{figure*} 
\centering 
\includegraphics[scale = 0.35]{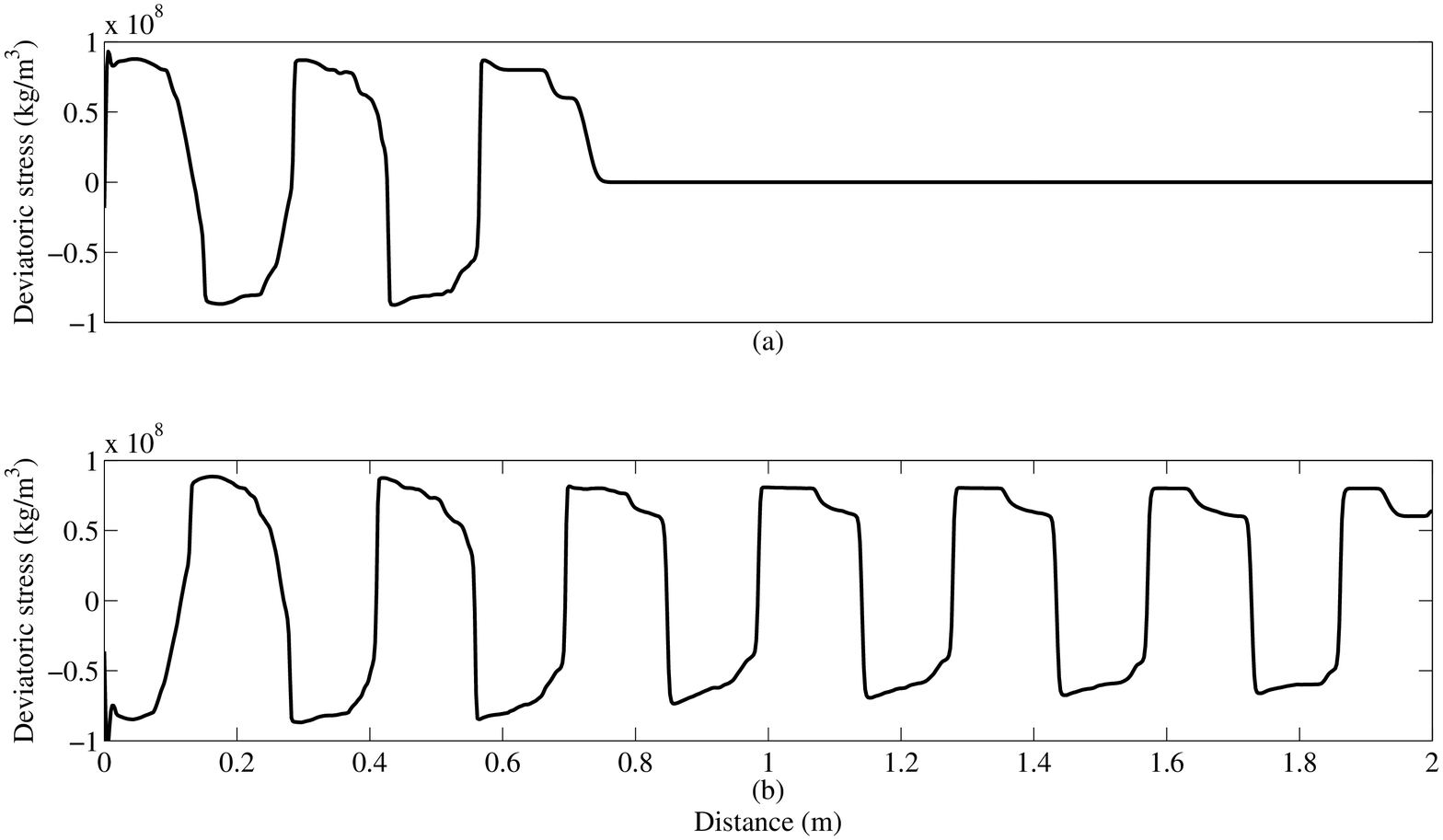}
\caption{Stress profiles of stress wave propagation in bulk elastic-hardening
plastic materials. The left end is applied with a sinusoidal force in 10 kHz
frequency. Amplitude of the applied force is larger than yielding stress of copper.}
\label{f:plastic9} 
\end{figure*}

\noindent
In this setting, we
consider a copper bar whose left end is forced by a sinusoidal ultrasonic force
which is 
\begin{align} 
F = F_0 \sin(\omega t), \label{e:np:100}
\end{align} 
with $F_0 = 200000 KN$ and $\omega = 10^5 rad/s$.
For simplicity, we set the initial conditions zero and the area
of the bar to unity. Since the external force is high enough to cause copper
yielding, elastic and plastic wave in the bar will be generated. It will propagate
to the other end. Stress wave and velocity wave profiles are shown in
Fig.(\ref{f:plastic8}) and Fig.(\ref{f:plastic9}). When elastic and plastic wave 
propagates to the right, wave shapes have been changed. Elastic
tensional and compressional waves and plastic tensional and compressional
waves take place during the same time. The difference of elastic wave speed and 
plastic wave speed makes the wave profile continually changing.
 It's also noticed that the original symmetrical input
wave profiles have been changed to unsymmetrical elasto-plastic wave profiles.

This asymmetry could be caused by the
asymmetrical effective stress-strain curve. That is,  in the copper constitutive model, yielding compression stress
is different from yielding tensile stress.
Because of the nonlinearity of the governing equations, the density and stress profiles shown in 
Fig.(\ref{f:plastic8}) and Fig.(\ref{f:plastic9}) no longer keep sinusoid shape. The
oscillations of many other frequency waves are added due to nonlinear effect. The final
wave shape is becoming square wave.

\section{Limitation and Conclusion}
\subsection{Limitation}
To understand stress wave propagation in multiple mediums and multiple phases such as solid-fluid interface and large amplitude elastic-plastic flow is a very challenge subject. The difficulties are really coming from several aspects. 
Firstly, governing equations in Eulerian frame consist of conservation of mass, momentum, energy, and many transport equations. They are very complicated nonlinear partial differential equations, which are all coupled together. This actually makes analytically solvable problems unsolvable. In vibration analysis or quantum mechanics, we are trying very hard to decouple multiple differential equations by either introducing eigenvectors and eigenvalues or assuming separation variables. The original purpose of solving conservation laws in Eulerian frame in fixed volume and fixed mesh is to deal with large deformation and plastic flow which may be challenging for Lagrangian method. But it introduces so many nonlinearities into the governing equations which makes problems even harder.
The other limitation of governing equations in continuum mechanics scale is that no micro-structure evolutions are considered. When we talk about conservation of momentum, we talk about linear momentum. However, in continuous scale, We don't know whose linear momentum this is. Does this linear momentum belongs to dislocations, electrons, or molecules? If this linear momentum belongs to a finite volume, what is the size of finite volume? What if there exist many cracks or many interfaces inside the finite volume. Linear momentum may not be conserved. This linear momentum can be changed to surface motion on the cracks or angular momentum of molecules. In large deformable mediums, there is a chance that linear momentum can be exchanged with angular momentum which is not considered here. When Isaac Newton introduced the second law, e.g., $F=ma$, he was working on particles. He never expected people will extend his theory to deformable materials. It is Cauchy, Lagrange, Navier, Stokes, and Euler who introduced particle mechanics to deformable materials which forms solid mechanics and fluid dynamics. When they built the theories, none of them considered micro-structures or quantum energy levels of solid deformable materials. They told us that mass, momentum, and energy are conserved. But they didn't tell us how fast they are conserved, where they will be conserved , when they will be conserved. Lots of concepts and theories still need be validated. Validation will be understood through how atoms and molecules bind themselves together. A path from micro-structure to macro-structure need be set up \citep{yang_digest_2023}. Quantum theory and statistical mechanics will be the answer.

Secondly, to build a large amplitude elastic-plastic constitutive model is very difficult. Traditional elastic plastic yielding criteria such as Von Mises or Drucker-Prager models are all rate independent plastic models. Their usage with conservation laws to model high speed impact is questionable. Their mathematical structures can hardly be cast into hyperbolic differential equation. In reality, when plastic deformation happens, dislocation will be created or annihilated. Dislocation density in crystalline structures \citep{yang_revisit_2020} or defect density in amorphous structures \citep{yang_mathematical_2019} will be changed and accumulated. In large deformation elastic and plastic $F_eF_p$ model, no micro-structures are considered. Writing $F_eF_p$ model into a transport equation coupled with conservation laws is almost impossible.
Thirdly, Discretization of large deformation elastic and plastic material model in numerical method such as CESE method is very complex. Some part of the constitutive model has to be treated as a source term and leave the other part in the hyperbolic form. With this treatment, time marching speed in numerical code will be different from wave speed in plastic zone. So if we have a very good numerical method such as CESE method, mathematical treatment of plastic deformation is not perfect. The numerical results can go to a different direction.
\subsection{Conclusion} 
\label{s:conclusion}
In this article, we pay more attention to understand plastic modeling instead of numerical method.
The isothermal hyperbolic model for stress wave in elastic-plastic solid does
not include energy conservation equation. Equation of state employed
relates pressure and density, without considering internal energy. We applied
the isothermal model to simulate low-speed impact problems. The numerical results are validated by comparing to an
analytical solution, which was derived by using a more comprehensive equation
state with the thermal effect.

 The above results show that the isothermal model developed in the present
paper correctly predicts the elastic-plastic wave propagation in low and moderate-speed impact. Thus, if
temperature is not of concern in a low-impact-speed problem, process which
might be able to assumed as a isothermal problem due to a slight temperature
change and a low material particle speed, one may use the isothermal model to
simulate process instead of complete model including the thermal effect.

\section*{Acknowledgement}
The help and generous support from Dr. John Sheng-tao Yu, Dr. Steven E Bechtel, and Dr. Minghao Tsai are greatly appreciated.
The authors also wish to acknowledge the generous support of this work by National Science Foundation Grant DMI-0600060.
\bibliographystyle{elsarticle-num}
\bibliography{fingon_dft}

\pagebreak
\begin{appendices}
\renewcommand{\theequation}
{A-\arabic{equation}}

{\bf {\em One dimension plastic strain wave analysis}}
\vspace{0.075in}

\begin{figure}
\centering
\includegraphics[scale = 0.5]{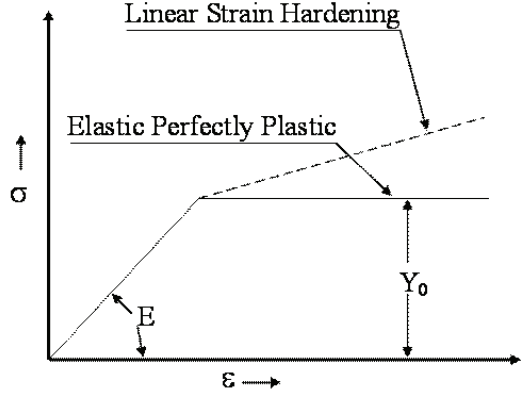}
\caption{uniaxial stress state curve for elastic, plastic material}
\label{f:plastic41}
\end{figure}

\noindent 
The conventional uniaxial stress-strain curve, as depicted by the
idealized models of Fig. (\ref{f:plastic41}), does not adequately represent
the state of stress and strain to which a material is subjected under shock
loading. 

If we consider the case of plane waves propagating through a material where
dimensions and constraints are such that the lateral strains are zero, the
stress and strain curve takes on a different form. This situation is commonly
referred to as uniaxial strain. 

To understand why these changes occur, consider the stress-strain relationship
for one-dimensional deformation. In the general case the three principal
strains can be divided into an elastic and a plastic part: 

\begin{equation}
\epsilon_1 = \epsilon_1^e +\epsilon_1^p, \\ \epsilon_2 = \epsilon_2^e +
\epsilon_2^p, \\ \epsilon_3 = \epsilon_3^e + \epsilon_3^p, \label{e:cons:A1}
\end{equation}  

where the superscripts $e$ and $p$ refer to elastic and
plastic, respectively, and  the subscripts are the three principal directions.

In one-dimensional deformation
\begin{equation}
\epsilon_2 = \epsilon_3 = 0, \\
\epsilon_2^p = - \epsilon_2^e, \\
\epsilon_3^p = - \epsilon_3^e, \label{e:cons:A2}
\end{equation}  

The plastic portion of the strain is taken to be incompressible, so that

\begin{equation}
\epsilon_1^p + \epsilon_2^p + \epsilon_3^p = 0 \label{e:cons:A3}
\end{equation}

which gives

\begin{equation}
\epsilon_1^p = -\epsilon_2^p - \epsilon_3^p = -2 \epsilon_2^p \label{e:cons:A4}
\end{equation}

since $\epsilon_2^p = \epsilon_3^p$ due to symmetry. From Eq.(\ref{e:cons:A2})
we have that

\begin{equation}
\epsilon_1^p = 2 \epsilon_2^e. \label{e:cons:A5}
\end{equation}

And the total strain $\epsilon_1$ may be written as

\begin{equation}
\epsilon_1 = \epsilon_1^e + \epsilon_1^p = \epsilon_1^e + 2 \epsilon_2^e. 
\label{e:cons:A6}
\end{equation}

The elastic strain in terms of the stresses and elastic constants is given by 

\begin{align}
\begin{aligned}
\epsilon_1^e &= \frac{\sigma_1}{E}-\frac{\nu}{E}(\sigma_2 + \sigma_3) 
= \frac{\sigma_1}{E} -\frac{2\nu}{E}\sigma_2, 
\\
\epsilon_2^e &= \frac{\sigma_2}{E}-\frac{\nu}{E}(\sigma_1 + \sigma_3) 
= \frac{(1-\nu)}{E}\sigma_2 -\frac{\nu}{E}\sigma_1, 
\\
\epsilon_3^e &= \frac{\sigma_3}{E}-\frac{\nu}{E}(\sigma_1 + \sigma_2) 
= \frac{1- \nu}{E}\sigma_3 -\frac{\nu}{E}\sigma_1.   
\label{e:cons:A71}
\end{aligned}
\end{align}
where $\sigma_2 = \sigma_3$ has been assumed due to symmetry. 
Combining Eq.(\ref{e:cons:A6}) and Eq.(\ref{e:cons:A71}),  we get

\begin{equation} 
\epsilon_1 = \frac{\sigma_1(1-2 \nu)}{E} + \frac{2
\sigma_2(1-2\nu)}{E}. 
\label{e:cons:A8} \end{equation}

And we know that the Von Mises plasticity yielding condition for this case is 

\begin{equation}
\sigma_1 - \sigma_2 = Y_0, \label{e:cons:A9}
\end{equation} 

where $Y_0$ is the static yield strength. If Eq.(\ref{e:cons:A8}) and
Eq.(\ref{e:cons:A9}) are combined to remove $\sigma_2$, we get

\begin{equation} 
\sigma_1 = \frac{E}{3(1-2\nu)}\epsilon_1 + \frac{2}{3}Y_0 = K
\epsilon_1 + \frac{2 Y_0}{3}, 
\label{e:cons:A10} 
\end{equation}

where $K = E/3(1-2\nu)$ is called the bulk modulus. Eq.(\ref{e:cons:A10}) is
valid only above yielding. The uniaxial-strain elastic stress-strain relation can
be easily got by considering first two equations of Eq.(\ref{e:cons:A71}) with
$\epsilon_2^e = 0$

\begin{equation} 
\sigma_1 = \frac{E(1-\nu)}{(1+\nu)(1-2\nu)} \epsilon_1 = (K +
\frac{4\mu}{3})\epsilon_1. 
\label{e:cons:A11} 
\end{equation}

\begin{figure}
\centering
\includegraphics[scale = 0.35]{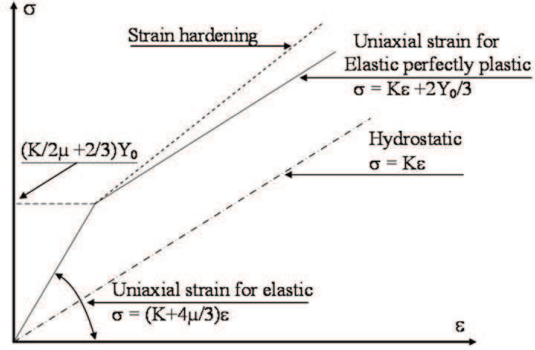}
\caption{uniaxial strain state curve for elastic, plastic material}
\label{f:plastic51}
\end{figure}

The uniaxial stress-strain curve for an elastic and plastic material is shown in
Fig.(\ref{f:plastic51})

Linear momentum equation and general one-dimensional stress-strain relation may be combined to be

\begin{equation} 
\frac{\d^2 u}{\d t^2} = \frac{1}{\rho_0}\frac{\d \sigma_1}{\d
x} = \frac{1}{\rho_0}\frac{\d \sigma_1}{\d \epsilon_1} \frac{\d^2 u}{\d x^2},
\label{e:cons:A12} 
\end{equation}

where relation $\epsilon_1 = \d u/\d x$ have been used.
So the wave speed of a strain increment is given by

\begin{equation} 
c = (\frac{1}{\rho_0} \frac{\d \sigma_1}{\d
\epsilon_1})^{0.5}. 
\label{e:cons:A13} 
\end{equation}

Combine Eq.(\ref{e:cons:A11}) and Eq.(\ref{e:cons:A13}), an elastic wave in
uniaxial strain configuration is

\begin{equation}
c_e = \sqrt{\frac{K + 4\mu/3}{\rho_0}}. \label{e:cons:A14}
\end{equation}

Combine Eq.(\ref{e:cons:A10}) and Eq.(\ref{e:cons:A13}), an plastic wave in
uniaxial strain configuration is

\begin{equation}
c_p = \sqrt{\frac{K}{\rho_0}}. \label{e:cons:A15}
\end{equation}

The wave speeds shown in Eq. (\ref{e:cons:A14}) and Eq.(\ref{e:cons:A15}) are
consistent with the eigenvalues shown in Eq. (\ref{e:cons:eigenvalue1}) where
$\beta = 0$ correspond to elastic wave speed and $\beta =1$,$B_{SH}= 0$
correspond to plastic wave speed.
\end{appendices}  

\end{document}